\documentclass[times]{nlaauth}
\usepackage{todonotes}
\usepackage{amsfonts}
\usepackage{amsmath}
\usepackage[colorlinks,bookmarksopen,bookmarksnumbered,citecolor=red,urlcolor=red]{hyperref}

\newcommand{\beq}{\begin{equation}}
\newcommand{\eeq}{\end{equation}}
\newcommand{\beqn}{\begin{equation*}}
\newcommand{\eeqn}{\end{equation*}}

\newcommand\BibTeX{{\rmfamily B\kern-.05em \textsc{i\kern-.025em b}\kern-.08em
T\kern-.1667em\lower.7ex\hbox{E}\kern-.125emX}}

\def\volumeyear{2013}

\begin{document}

\title{A new level-dependent coarse grid correction scheme for \\ indefinite Helmholtz problems}

\author{Siegfried Cools, Bram Reps, Wim Vanroose$^*$}

\address{Department of Mathematics and Computer Science, University of Antwerp, Middelheimlaan 1, 2020 Antwerp, Belgium \\ Intel ExaScience Lab, Kapeldreef 75, B-3001 Leuven, Belgium (B. Reps). \vspace{-0.8cm}}

\corraddr{siegfried.cools@ua.ac.be, wim.vanroose@ua.ac.be}


\begin{abstract}
In this paper we construct and analyze a level-dependent coarse grid correction scheme for indefinite Helmholtz problems. This adapted multigrid method is capable of solving the Helmholtz equation on the finest grid using a series of multigrid cycles with a grid-dependent complex shift, leading to a stable correction scheme on all levels. It is rigorously shown that the adaptation of the complex shift throughout the multigrid cycle maintains the functionality of the two-grid correction scheme, as no smooth modes are amplified in or added to the error. In addition, a sufficiently smoothing relaxation scheme should be applied to ensure damping of the oscillatory error components. Numerical experiments on various benchmark problems show the method to be competitive with or even outperform the current state-of-the-art multigrid-preconditioned Krylov methods, like e.g.~complex shifted Laplacian (CSL) preconditioned GMRES.
\end{abstract}


\maketitle

\vspace{-0.3cm}

\section{Introduction}

Originally introduced as a theoretical tool by Fedorenko in 1964 \cite{fedorenko1964speed} and later adopted as a solution method by Brandt in 1977 \cite{brandt1977multi}, the multigrid method is known to be a particularly fast and scalable solver for the large systems of equations arising from the discretization of multi-dimensional Poisson or positive definite Helmholtz equations, see \cite{brandt1986multigrid,briggs2000multigrid,stüben1982multigrid,trottenberg2001multigrid}. Using a negative shift, however, the Helmholtz discretization matrix becomes distinctly indefinite causing multigrid convergence to deteriorate. This break-down is due to near-to-zero eigenvalues in the operator spectrum on some intermediate level in the multigrid hierarchy, which destroy the functionality of the standard correction scheme, as shown by Elman, Ernst \& O'Leary in \cite{elman2002multigrid} and Ernst \& Gander in \cite{ernst2010difficult}.

Motivated primarily by geophysical applications \cite{riyanti2006new}, rising interest in the development of fast solvers for the indefinite Helmholtz equation over the past decade has yielded a broad range of solution methods for this non-trivial problem, an overview of which can be found in \cite{ernst2010difficult}. Due to their applicability on a variety of problems, preconditioned Krylov subspace methods are currently among the most popular Helmholtz solvers. Here the Krylov subspace method functions as an outer iteration and a direct (ILU) or iterative (multigrid) method is used to (approximately) solve the preconditioning system in every step. Being crucial to the performance of the governing Krylov method, significant research has been performed over the construction of a good preconditioner. Past and recent work includes the wave-ray approach \cite{brandt1997wave}, the idea of separation of variables \cite{plessix2003separation}, algebraic multilevel methods \cite{bollhoefer2009algebraic} and a transformation of the Helmholtz equation into an advection-diffusion-reaction problem \cite{haber2011fast}. The key ideas of the current paper show some resemblances to the work done in \cite{maclachlan2008algebraic}, where the use of complex valued preconditioning operators are analyzed in an AMG setting.

Another multigrid-based preconditioner was suggested in \cite{elman2002multigrid}, which has led to the development of the so-called shifted Laplacian preconditioners. The first papers on these preconditioners are \cite{bayliss1983iterative} and \cite{lairdpreconditioned}, in which a Laplace operator with a real shift was proposed as a preconditioner instead of the approximate inverse of the original Helmholtz operator. However, still being unable to efficiently solve very high wavenumber problems using this approach, the idea was reinvented and extended to the Complex Shifted Laplacian (CSL) by Erlangga, Vuik and Oosterlee in \cite{erlangga2004class,erlangga2006novel} and analyzed further in \cite{erlangga2006comparison,vangijzen2007spectral}. Leading to satisfactory convergence and scalability results on highly indefinite problems, the complex shifting of the original problem operator was generalized in \cite{erlangga2008multilevel,erlangga2009algebraic,airaksinen2009algebraic,osei2010preconditioning}.

In \cite{reps2010indefinite} it was shown that scaling the wavenumber by a complex value is equivalent to scaling the grid distance, thereby rotating the spectrum around its most negative real point instead of translating it up or down in the complex plane. The effects of this Complex Stretched Grid (CSG) transformation on MG-preconditioned Krylov solvers were analyzed to some extent in \cite{reps2011analyzing}. The resulting preconditioner is a Helmholtz operator discretized on a complex-valued grid, which is proven to be particularly efficient when complex-valued grid distances are already used to implement advanced absorbing boundary conditions like Exterior Complex Scaling (ECS) \cite{aguilar1971class} or Perfectly Matched Layers (PML) \cite{berenger1994perfectly}.

In this paper we aim to solve the indefinite Helmholtz equation using a new level-dependent multigrid scheme. The method adopts the key idea from CSL/CSG, which is plugged into the multigrid correction scheme, gradually shifting the eigenvalues away from the origin in the coarsening process. This results in a multigrid scheme which, contrarily to the above methods, can be applied \emph{directly} as a solver to the Helmholtz equation, instead of merely suiting preconditioning purposes. Numerical results show the method to be competitive with (or even outperform) the current state-of-the-art MG-preconditioned Krylov solvers.

To complement this adapted correction scheme, a stable relaxation scheme with strong smoothing properties should be used to effectively damp the oscillatory error components. It is well-known that for Helmholtz problems standard smoothers such as weighted-Jacobi or Gauss-Seidel do not necessarily resolve high-frequency modes and might even display divergent behaviour. Considering this observation we recommend using the more advanced GMRES(3) method as a smoother substitute. First suggested as a replacement for standard multigrid smoothers in \cite{elman2002multigrid}, it was recently shown in \cite{calandra2011two} and \cite{reps2011multigrid} that replacing the standard smoother by GMRES(3) yields very satisfactory multigrid convergence results when solving Helmholtz equations.

The paper is organized as follows. In Section 2 we briefly recall the properties of both the Complex Shifted Laplacian (CSL) and Complex Stretched Grid (CSG) preconditioners. Additionally, we show the stability and functionality of these schemes using the spectral analysis from \cite{elman2002multigrid}. In Section 3 we then define the new level-dependent correction scheme, in which we differentiate between a CSL- and CSG-based variant. The level-dependent scheme is analyzed through a standard two-grid spectral analysis, emphasizing the complementary action of the correction scheme and the smoother. Numerical results supporting the theory are shown in Section 4, where the new level-dependent scheme is extensively tested on a variety of benchmark problems, ranging from the academic constant wavenumber model problem to more hard-to-solve quantum mechanical models featuring evanescent waves. Finally, alongside a discussion on the subject, conclusions are drawn in Section 5.

\section{Overview of complex shifted preconditioners}

\noindent It is the aim of this work to effectively solve the discretized Helmholtz system
\beq \label{eq:HH}
(-\Delta^h - k^2) \, u^h = \chi^h, \qquad \text{on~} \Omega^h \subset \mathbb{R}^d,
\eeq
with dimension $d \geq 1$, where $\Delta$ is the Laplace operator, $\chi$ represents the right-hand side or \emph{source term}, and $k \in \mathbb{R}$ is known as the (possibly spatially dependent) \emph{wavenumber} that might cause the system to become highly indefinite. Note that for the analysis we generally restrain ourselves to the one-dimensional problem formulation with Dirichlet boundary conditions and a constant wavenumber $k$, as is common practice in Helmholtz literature. The theoretical results obtained in this work can, however, readily be extended to more general problem settings in higher dimensions. Successful experiments in higher spatial dimensions (see Section 4) will provide a solid foundation for this statement. For completeness, we note that the subsequent analysis is essentially based upon a geometric multigrid setting, using rediscretization to define the coarse grid operators.

\subsection{Two-grid spectral analysis}
As a multigrid cycle is essentially  a recursive embedding of two-grid correction schemes on consecutive coarser grids, flanked on both sides by a pre- and postsmoothing operator, an analysis of the two-grid scheme is crucial in understanding the full action of the multigrid method. The two-grid correction scheme for an error $e^h \in \Omega^h$ is given by the matrix operator (see \cite{briggs2000multigrid}, \cite{stüben1982multigrid})
\beq \label{eq:standardtg}
TG = \left(I - I_{2h}^h (A^{2h})^{-1} I_h^{2h} A^h\right) 
\eeq
where $I$ is the identity matrix, $A^h$ and $A^{2h}$ are the fine and coarse grid discretization matrices respectively, $I_{2h}^h$ is the interpolation or prolongation operator, and $I_h^{2h}$ is the restriction operator. The standard intergrid operators used in this text are (bi-/tri-)linear interpolation and full weighting restriction. Note that any vector $v^h \in \Omega^h$, specifically the error $e^h$, can be written as a linear combination of the eigenvectors of $A^h$
\beq
e^h = \sum_{l = 1}^{n} a_l w_l^h,
\eeq
with $a_l \in \mathbb{R}$, where it is known that the eigenvectors of $A^h$ are given componentwise by
\beq
w_{l,j}^h = \sin\left(\frac{lj\pi}{n}\right), \quad 1 \leq j \leq n-1.
\eeq
Hence it suffices to study the action of the coarse grid scheme on the eigenvectors $w_l^h$. It is common for the analysis to distinguish between smooth eigenmodes $w_l^h$ with $1 \leq l < n/2$ and the corresponding oscillatory eigenmodes $w_{l'}^h$, where $l'= n - l$. Defining the constants $c_l = \cos^2(\frac{l\pi}{2n})$ and $s_l = \sin^2(\frac{l\pi}{2n})$, the following general expressions hold
\beq \label{eq:tgsmoo}
TG \, w_l^h = w_l^h - c_l \frac{\lambda_l^h}{\lambda_l^{2h}}(c_l w_l^h - s_l w_{l'}^h),
\eeq
and
\beq
TG \, w_{l'}^h = w_{l'}^h - s_l \frac{\lambda_{l'}^h}{\lambda_l^{2h}}(c_l w_l^h - s_l w_{l'}^h),
\eeq
where $\lambda_l^h$ and $\lambda_{l'}^h$ are the eigenvalues corresponding to $w_l^h$ and $w_{l'}^h$ respectively. The classical two-grid analysis of Ernst, Elman and O'Leary \cite{elman2002multigrid} follows as a smooth mode limit case from expression (\ref{eq:tgsmoo}), as for $l \ll n/2$ we have $c_l \approx 1$ and $s_l \approx 0$ implying
\beq \label{eq:olea}
TG \, w_l^h \approx \left(1 - \frac{\lambda_l^h}{\lambda_l^{2h}}\right) w_l^h.
\eeq
For notational simplicity we now denote $\lambda_l^h = \lambda^h$ and $\lambda_l^{2h} = \lambda^H$. For the indefinite Helmholtz equation it is well-known that $\lambda^h = \lambda^h_L - k^2$, where $\lambda^h_L \in [0,2^{d+1}/h^2]$ is the corresponding eigenvalue of the discretized $d$-dimensional Laplace operator $L = -\Delta^h$. As discussed in \cite{elman2002multigrid} a quasi-zero denominator in expression (\ref{eq:olea}) due to $\lambda^{H}_L \approx k^2$ and/or oppositely signed eigenvalues $\lambda^h$ and $\lambda^H$ can lead to two-grid instability, i.e.
\beq
\left|1-\frac{\lambda^h}{\lambda^H}\right| = \left|1-\frac{\lambda_L^h - k^2}{\lambda_L^H - k^2}\right| \gg 1.
\eeq
Indeed, this observation acted as the main motivation for the development of the CSL and CSG preconditioners, see further. Note, however, that the eigenvalues $\lambda_L^h$ and $\lambda_L^H$ corresponding to the smoothest eigenmodes $w_l^h$ resp.~$w_l^{2h}$ with $l \ll n/2$ are relatively close to zero in comparison to $k^2$. Indeed, we have that both $\lambda_L^h \ll k^2$ and $\lambda_L^H \ll k^2$, as $l \ll n/2$ implies that $s_l \approx 0$ and $c_l \approx 1$ such that
\begin{align}
\lambda_L^h &= \frac{4d}{h^2} \sin^2\left(\frac{l\pi}{2n}\right) =  \frac{4d}{h^2} s_l \approx 0 \ll k^2, \\
\lambda_L^H &= \frac{d}{h^2} \sin^2\left(\frac{l\pi}{n}\right)  =  \frac{4d}{h^2} \sin^2\left(\frac{l\pi}{2n}\right) \cos^2\left(\frac{l\pi}{2n}\right) = \frac{4d}{h^2} s_l c_l \approx 0 \ll k^2.
\end{align}
Hence the smoothest eigenmodes are always mapped onto zero by the $TG$ operator, as (\ref{eq:olea}) implies that
\beq
TG \, w_l^h \approx  \left(1-\frac{\lambda_L^h - k^2}{\lambda_L^H - k^2}\right) w_l^h = \left(1-\frac{\lambda_L^h/k^2 - 1}{\lambda_L^H/k^2 - 1}\right) w_l^h \approx 0 \cdot w_l^h.
\eeq
Note that this mapping of the smoothest modes onto zero is a desirable property for any coarse grid correction scheme, as it guarantees optimal collaboration between correction scheme and smoother.

\subsection{Complex shifted Laplacian}
In the 2004 paper \cite{erlangga2004class} Erlangga, Vuik and Oosterlee introduced a solution to the occurrence of small-valued coarse grid eigenvalues that destroy two-grid stability by means of introducing a real and/or complex shift to the problem. The slightly perturbed Helmholtz problem 
\beq \label{eq:csl}
(-\Delta^h - (\alpha + \beta i)k^2) u^h = \chi^h
\eeq
proposed within \cite{erlangga2004class} can be solved efficiently using multigrid given a sufficiently large complex shift part $\beta$. A recent discussion on lower bounds for the shift parameter $\beta$ can be found in \cite{cools2012local}. Note that alternatively, the above expression can be reformulated in a slightly less convenient form as $(-\Delta^h - re^{i\theta}k^2) u^h = \chi^h$ such that the eigenvalues of the perturbed discretization matrix $\tilde{A}^h$ are given by $\lambda^h = \lambda_L^h - re^{i\theta}k^2$. The two-grid correction scheme on this perturbed problem is
\beq 
TG = \left(I - I_{2h}^h (\tilde{A}^{2h})^{-1} I_h^{2h} \tilde{A}^h\right),
\eeq
where $\tilde{A}^h = L - re^{i\theta}k^2$. Note that the eigenvectors $w_l^h$ of the original Helmholtz problem (\ref{eq:HH}) are preserved by shifting the Laplace operator. Hence, following expression (\ref{eq:olea}) for the smooth modes, the denominator of the fraction 
\beq \label{eq:csltgw}
\left|1-\frac{\lambda_{L}^h-re^{i\theta}k^2}{\lambda_{L}^H-re^{i\theta}k^2}\right|
\eeq
will never approach zero as $|\lambda^H| = |\lambda_L^H - re^{i\theta}k^2| > |\mathcal{I}(re^{i\theta} k^2)| > 0$ for all coarse grid eigenvalues $\lambda_L^H$ due to the complex part. Similarly to the unperturbed problem, the smoothest eigenvalues are mapped onto zero by the $TG$ operator, as for $l \ll n/2$ we have $\lambda_L^h \ll k^2$ and $\lambda_L^H \ll k^2$ such that
\beq \label{eq:csltg}
TG \, w_l^h =  \left(1-\frac{\lambda_L^h - re^{i\theta}k^2}{\lambda_L^H - re^{i\theta}k^2}\right) w_l^h = \left(1-\frac{\lambda_L^h/k^2- re^{i\theta}}{\lambda_L^H/k^2 - re^{i\theta}}\right) w_l^h \approx 0 \cdot w_l^h.
\eeq
Thus for a sufficiently large complex shift, the two-grid correction scheme is stable for the shifted problem (\ref{eq:csltgw}) and maps the smooth eigenvalues near the origin as requested (\ref{eq:csltg}). However, despite its ability to efficiently solve the perturbed problem (\ref{eq:csl}), the original Helmholtz problem (\ref{eq:HH}) remains unsolvable by multigrid, and the multigrid solution method for the CSL problem (\ref{eq:csl}) can only be used as a preconditioner solver to the original Helmholtz problem. The remaining preconditioned matrix-vector system is most commonly solved by means of a Krylov subspace method like e.g.~restarted GMRES \cite{saad1986gmres}.

\begin{figure}[t]
\begin{center}
\includegraphics[width=6cm]{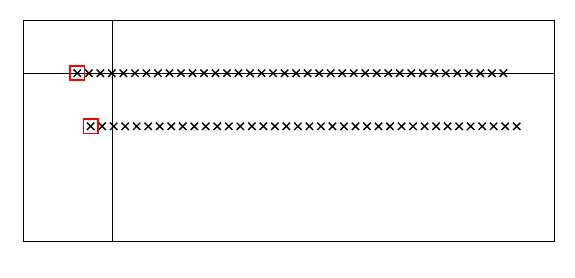}
\includegraphics[width=6cm]{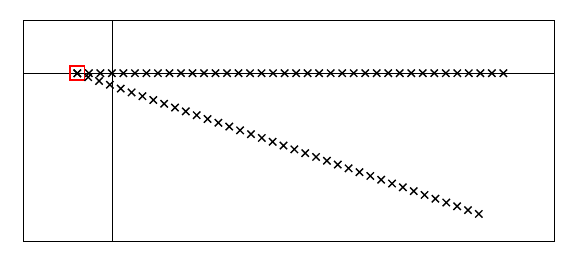}
\vspace{-0.2cm}
\caption{Schematic representation of the spectra of the Dirichlet bounded discretization matrix operator $A^{h}$ and the perturbed CSL (left) and CSG (right) operators $\tilde{A}^{h}$. Note how the CSG perturbation does not change the location of the leftmost (smoothest) eigenvalue, whereas CSL translates the entire spectrum.}
\label{fig:specsketch}
\end{center}
\end{figure}

\subsection{Complex stretched grid}
Introduced as an alternative to Complex Shifted Laplacian in \cite{reps2010indefinite}, the Complex Stretched Grid (CSG) preconditioner can be seen as an equivalent method of perturbing the Helmholtz problem to resolve the two-grid correction instability issue. As opposed to CSL, however, the existing shift $-k^2$ is left unchanged, whereas the discretization grid is rotated into the complex plane
\beq \label{eq:csghelm}
(-\Delta^h e^{-i\theta} - k^2) u^h = \chi^h.
\eeq
Note that we have deliberately chosen $r = 1$, inducing a pure rotation of the grid around $-k^2 \in \mathbb{R}$ (see Figure \ref{fig:specsketch}), such that the distance between two grid points is left unchanged $|\tilde{h}| = |he^{-i\frac{\theta}{2}}| = |h|$. To stress the similarity between CSL and CSG, note that the above CSG problem formulation (\ref{eq:csghelm}) is equivalent to the CSL problem
\beq \label{eq:csghelmprime}
(-\Delta^h - k^2 e^{i\theta}) u^h = \chi^h e^{i\theta}.
\eeq
As discussed in \cite{reps2010indefinite} this equivalence between CSL and CSG generally holds. The CSL wavenumber perturbation factor $1 + \beta i$, which induces no additional real shift and is therefore commonly used in Helmholtz literature, see e.g.~\cite{erlangga2006novel}, can instantly be transferred to a CSG setting by choosing $r = 1/\cos\theta$ and vice versa. As a general remark, note that for small rotation angles $\theta \approx 0$, this implies $r \approx 1$. In the following it is assumed, without loss of generality, that $r \equiv 1$.

The eigenvalues of the resulting perturbed CSG problem are given by $\lambda^h = \lambda_L^h e^{-i\theta} - k^2$. Referring once more to (\ref{eq:olea}), it follows from an argument similar to (\ref{eq:csltgw}) that the denominator in 
\beq 
\left|1-\frac{\lambda_{L}^he^{-i\theta}-k^2}{\lambda_{L}^He^{-i\theta}-k^2}\right|
\eeq
cannot approach zero when the rotation angle $\theta$ is sufficiently large, hence solving the two-grid instability issue. Additionally, the smooth mode eigenvalues $\lambda_l^h$ with $l \ll n/2$ are mapped onto zero by the two-grid operator as can be readily derived from (\ref{eq:olea}) in analogy to (\ref{eq:csltg}). This implies that the perturbed CSG problem can also be solved very efficiently using multigrid. Analogously to the CSL method, however, the Complex Stretched Grid approach can \emph{only} be used as a preconditioner to a general Krylov method, as the original Helmholtz problem cannot be solved by solely using multigrid.

\subsection{Influence of boundary conditions on the spectral analysis}

\begin{figure}[t]
\begin{center}
\includegraphics[width=6cm]{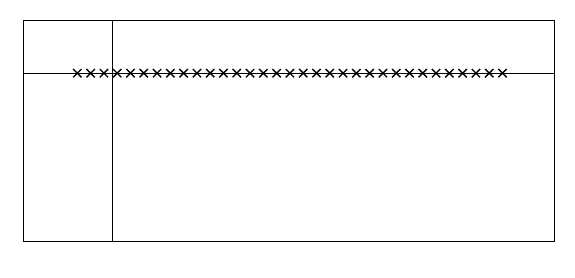}
\includegraphics[width=6cm]{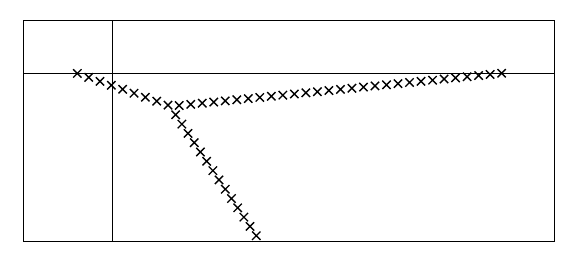}
\vspace{-0.2cm}
\caption{Schematic representation of the spectra of the Dirichlet (left) and ECS (right) bounded discretization matrix operator $A^{h}$. The use of ECS boundary conditions implies a natural damping, which is apparent from the spectrum.}
\label{fig:specsdirecs}
\end{center}
\end{figure}

We note that for convenience the above (and following) analysis assumes the use of Dirichlet-type conditions at the numerical boundaries of the problem. Throughout the Helmholtz literature Dirichlet boundary conditions are often suggested to imply a `worst case scenario' for efficient iterative solution, as no natural damping occurs. As illustrated by Figure \ref{fig:specsdirecs}, the use of more advanced absorbing boundary conditions like e.g.~Exterior Complex Scaling (ECS) generally results in a more amenable spectrum. Given a sufficiently fine discretization, the eigenvalues of the purely Dirichlet bounded problem plausibly approach the origin on every level in the multigrid hierarchy, as all eigenvalues are located along the real axis. The use of ECS boundary conditions, however, implies a natural shift of the eigenvalues into the negative half of the complex plane. This results in only a \emph{single} problematic (or `critical') level in the multigrid hierarchy, which corresponds to the eigenvalue(s) located near $4/h^2-k^2$ on the finest grid approaching (or crossing) the origin after multiple levels of coarsening (cfr. \cite{reps2011analyzing}). Note that the spectral analysis proposed above is mainly based upon the smoothest and most oscillatory eigenvalues, as representatives of all smooth/oscillatory eigenmodes respectively. Consequently, the spectral analysis remains generally valid for problems with absorbing boundary layers, as the leftmost (very smooth) and rightmost (highly oscillatory) eigenvalues of the ECS spectrum are intrinsically identical to their Dirichlet-bounded counterparts. 

\section{Level-dependent coarse grid correction scheme}

In this section we introduce a level-dependent correction scheme, based on the perturbation idea introduced by the CSL and CSG preconditioners. As opposed to the latter methods, however, the problem is now gradually perturbed throughout the grid hierarchy, ensuring both stability of the two-grid scheme on all levels as well as solution of the original problem on the finest grid. This leads to an adapted multigrid correction scheme that is \emph{directly} applicable as a solver for the original Helmholtz problem (\ref{eq:HH}). Additionally, it is shown that the small error introduced by gradually perturbing the problem does not undermine the functionality of the multigrid scheme.

\subsection{Definition and analysis}
We introduce the notion of a level-dependent correction scheme based on the CSL and CSG schemes. Assuming the total number of levels in the multigrid hierarchy equals $p \in \mathbb{N}_{0}\backslash\{1\}$, the two-grid scheme on the $m$-th level in the level-dependent multigrid hierarchy is defined as
\beq
TG^m = \left(I - I_{2h}^h (\tilde{A}^{2h})^{-1} I_h^{2h} \tilde{A}^h\right), \qquad 1 \leq m < p,
\eeq
where for the scheme based on CSL perturbation we state that
\begin{eqnarray} \label{eq:cslvariant}
\tilde{A}^{h} &=& -\Delta^h - k^2 e^{i\theta_{m-1}}, \notag \\
\tilde{A}^{2h} &=& -\Delta^{2h} - k^2 e^{i\theta_m},  \qquad 1 \leq m < p,
\end{eqnarray} 
or equivalently for the CSG-based variant
\begin{eqnarray} \label{eq:csgvariant}
\tilde{A}^{h} &=& -\Delta^h e^{-i\theta_{m-1}} - k^2, \notag \\
\tilde{A}^{2h} &=& -\Delta^{2h} e^{-i\theta_m} - k^2,  \qquad 1 \leq m < p.
\end{eqnarray} 
Here $h = 2^{m-1} h_f$, where $h_f$ is the grid distance on the finest grid. The level-dependent parameter $\theta_m$ is most commonly defined as $\theta_m = md\theta$ for a \emph{small} perturbation parameter $d\theta \in [0,\pi]$, however, other perturbation schemes are possible. In this work the per-level perturbation $d\theta$ is defined as $d\theta = \theta_{max}/p$ where the fixed angle $\theta_{max}$ is the maximal rotation angle of the coarsest level, usually chosen $\theta_{max} = \pi/6$. Note that, by definition, on the finest level (where $m=1$) we have $\tilde{A}^h = A^h = -\Delta^h-k^2$. This observation acts as our primary motivation for the adaptation of the two-grid scheme into a level-dependent shift or rotation angle rather than a fixed perturbation of the problem on all levels. The hierarchy of two-grid schemes proposed above is designed to ultimately solve the \emph{original} Helmholtz problem on the finest grid. 

On every grid the current solution $v^h$ is corrected by $v^h \leftarrow v^h + \tilde{e}^h$ using the interpolated solution $\tilde{e}^{h}$ of a slightly perturbed coarse grid problem as an approximation to the exact error. Note that even in a regular two-grid scheme, the transfer of the residual to its discrete coarse grid representation results in an \emph{approximation} of the fine grid problem (due to restriction and interpolation). In the level-dependent scheme, we minorly alter the problem definition on the coarse grid by adding a small shift or domain rotation, obtaining a slightly different approximation to the error. In doing so, we introduce a small additional error on the correction due to the difference in problem definition between coarse and fine grid. In the next paragraphs, however, it will be shown that the resulting two-grid scheme does not excite smooth modes in the correction term. Indeed, we will show that the additional error component introduced mainly consists of oscillatory eigenmodes, which are subsequently damped by the smoother, yielding an accurate overall error correction scheme.

\subsubsection{Two-grid analysis of the level-dependent scheme.} In the following paragraph we perform a standard two-grid analysis of the new level-dependent scheme, cfr.~Section 2. We distinguish between the CSL- and CSG-based level-dependent schemes as presented above. Focussing first on the CSL scheme and using expression (\ref{eq:olea}), we observe that the denominator of the fraction 
\beq \label{eq:tglvldep}
\left|1-\frac{\lambda_{L}^h-e^{i\theta_{m-1}}k^2}{\lambda_{L}^H-e^{i\theta_m}k^2}\right|
\eeq
never approaches zero for $\theta_m$ sufficiently large, resolving the instability issues of the standard two-grid scheme. Additionally, the amplification factor for the smoothest eigenmodes $w_l^h$ with $l \ll n/2$ can be calculated as follows
\begin{align} \label{eq:cslexp}
TG^m \, w_l^h &= \left(1-\frac{\lambda_{L}^h-e^{i\theta_{m-1}}k^2}{\lambda_{L}^H-e^{i\theta_m}k^2}\right)w_l^h = \left(1-\frac{\lambda_{L}^h/k^2-e^{i\theta_{m-1}}}{\lambda_{L}^H/k^2-e^{i\theta_m}}\right)w_l^h \notag \\
						&\approx (1-e^{-id\theta}) \, w_l^h ,
\end{align}
where we again assume that both $\lambda_{L}^h \ll k^2$ and $\lambda_{L}^H \ll k^2$ for the smoothest eigenmodes. Notably, the amplification factor of the smoothest eigenmodes depends on the parameter $d\theta$, which in general is very small, implying $(1-e^{-id\theta}) \approx 0$ such that the smoothest eigenvalue is indeed mapped closely to zero. From (\ref{eq:cslexp}) it can be derived that stability for the smooth eigenmodes is guaranteed as long as $d\theta < \pi/3$, which in practice will always be satisfied. Contrary to the standard two-grid schemes discussed in Section 2, however, the smoothest modes are not \emph{exactly} mapped onto zero by the CSL-based level-dependent $TG^m$ operator. Instead, the minimal distance between the origin and a smooth mode's projection under the CSL variant of $TG^m$ is $|1-e^{-id\theta}|$. 

Subsequently turning to the level-dependent CSG scheme, it is clear that this scheme is also generally stable as the denominator of 
\beq
\left|1-\frac{\lambda_{L}^he^{-i\theta_{m-1}}-k^2}{\lambda_{L}^He^{-i\theta_m}-k^2}\right|
\eeq
does not approach zero for $\theta_m > 0$ sufficiently large, due to the complex part. Assuming once more that the eigenvalues corresponding to the smoothest eigenmodes $w_l^h$ (with $l \ll n/2$) are relatively close to zero in comparison to $k^2$, such that $\lambda_{L}^h \ll k^2$ and $\lambda_{L}^H \ll k^2$, the action of $TG^m$ on a very smooth eigenmode is given by
\beq
TG^m \, w_l^h = \left(1-\frac{\lambda_{L}^he^{-i\theta_{m-1}}-k^2}{\lambda_{L}^He^{-i\theta_m}-k^2}\right) w_l^h = \left(1-\frac{(\lambda_{L}^h/k^2)e^{-i\theta_{m-1}}-1}{(\lambda_{L}^H/k^2)e^{-i\theta_m}-1}\right) w_l^h \approx 0 \cdot w_l^h.
\eeq
Consequently, one observes that using the CSG scheme the smoothest eigenvalues are mapped onto zero as requested, independently of the parameter $d\theta$. Note that this was not the case for the CSL variant of $TG^m$, see (\ref{eq:cslexp}). This subtle difference between the CSL and CSG level-dependent schemes can be clarified further by studying their fundamental operation on the spectrum of the discretization matrix $A^h$, see Figure \ref{fig:specsketch}: while CSL induces a shift of the entire fine grid spectrum over a distance defined by $d\theta$, causing the fine- and coarse grid eigenvalues in the level-dependent scheme to be distinctly different, CSG rotates the spectrum over an angle $-d\theta$ around the point $-k^2 \in \mathbb{R}$, such that the smoothest corresponding eigenvalues on the fine and coarse grid are nearly left unchanged. This signifies that for smooth eigenmodes the CSG level-dependent scheme (\ref{eq:csgvariant}) resembles the action of the standard coarse grid correction scheme (\ref{eq:standardtg}) somewhat closer than the CSL variant. In light of this observation, our preference in solving practical problems goes out to the CSG variant of the level-dependent scheme.

\subsubsection{Spectral analysis of the error components.} 
\begin{figure}[t!] 
\begin{center}
\includegraphics[width=6cm]{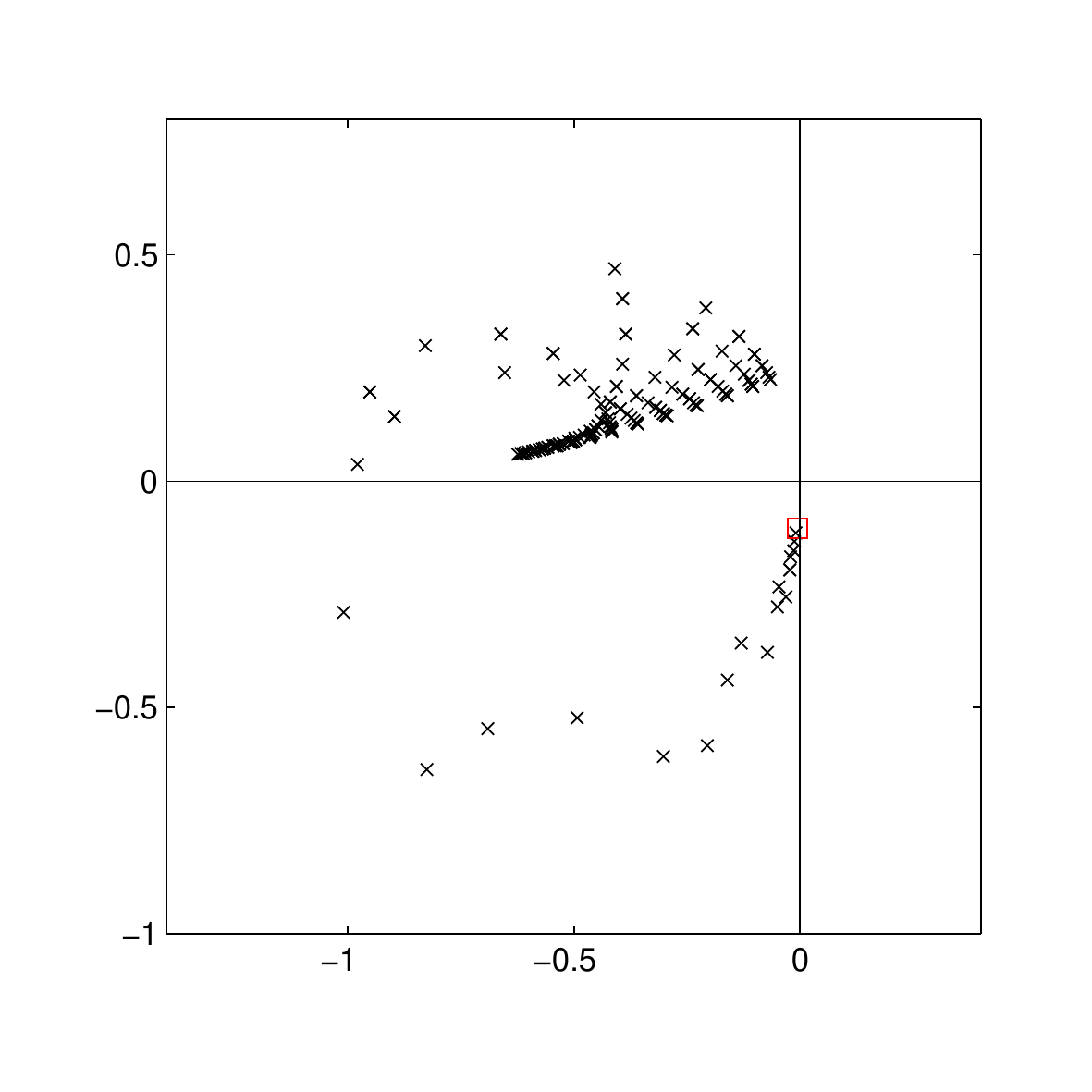}
\includegraphics[width=6cm]{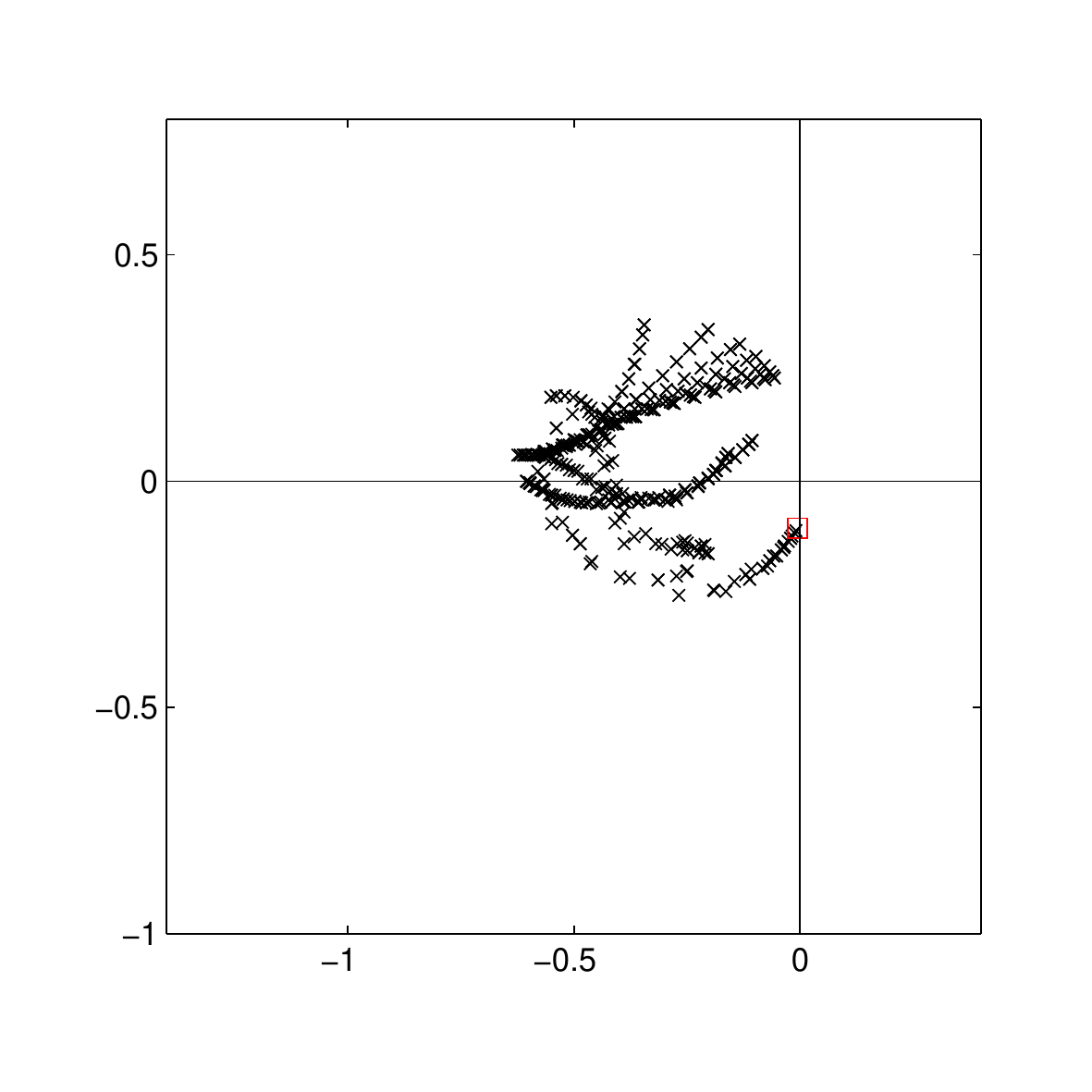}\\ 
\includegraphics[width=6cm]{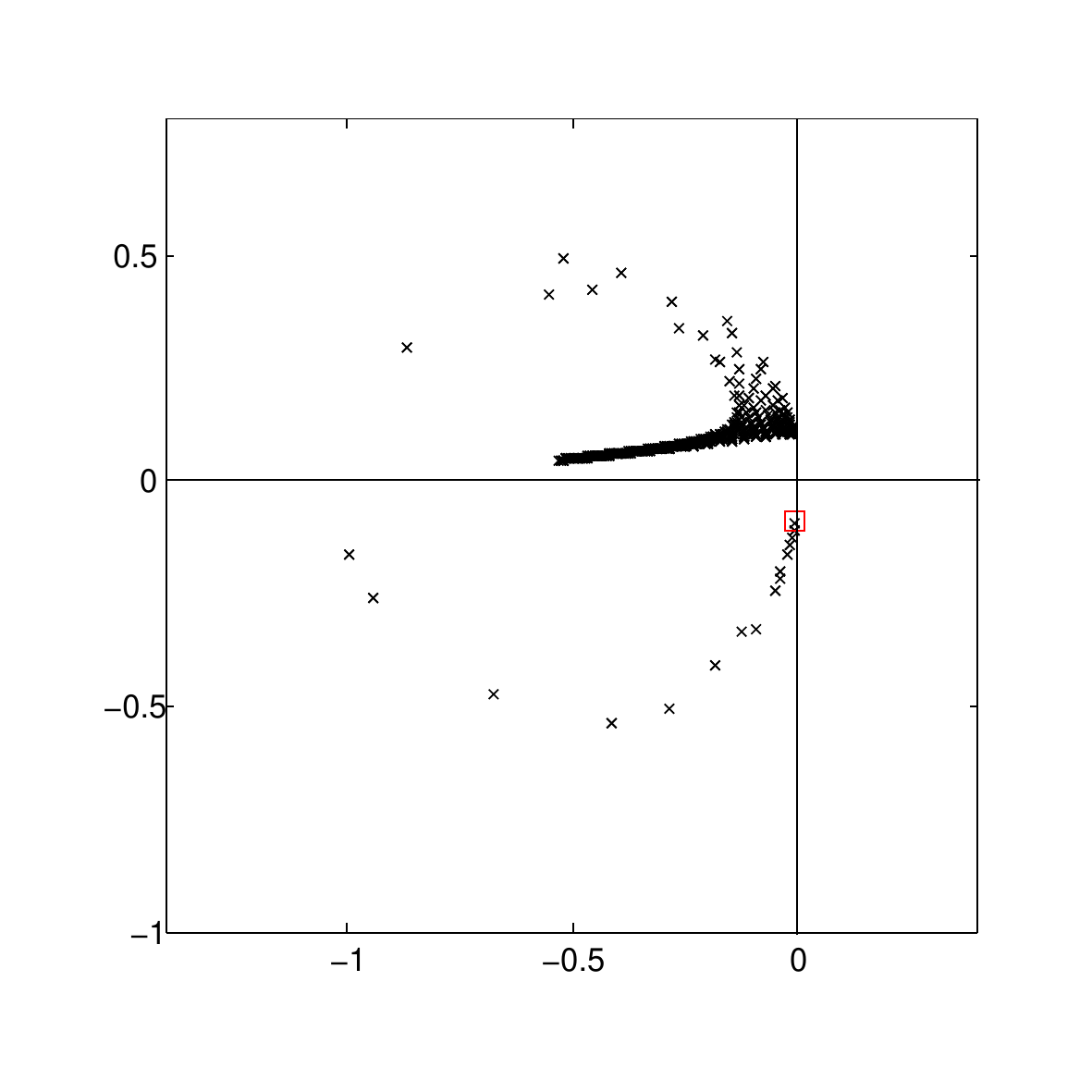} 
\includegraphics[width=6cm]{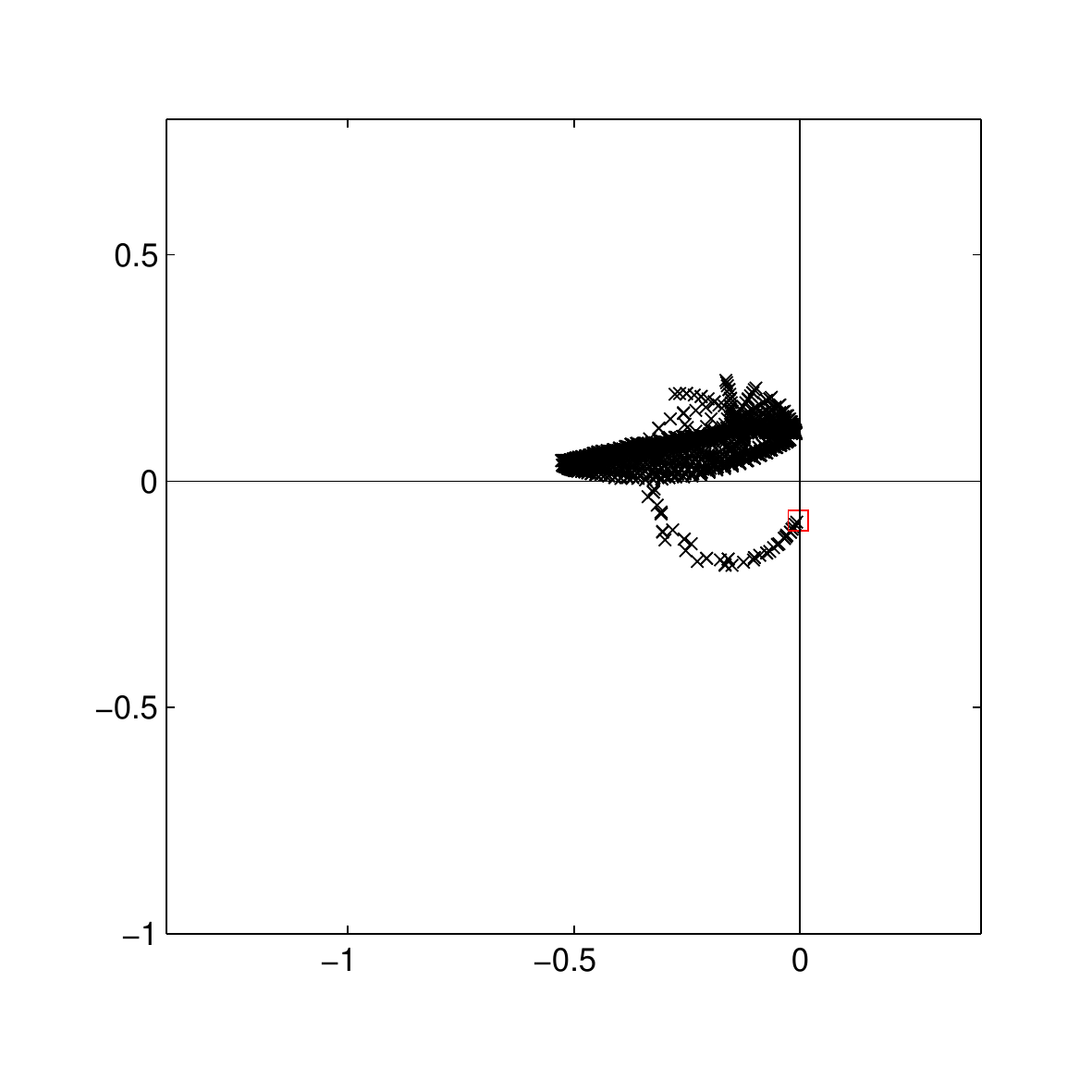}\\ 	
\includegraphics[width=6cm]{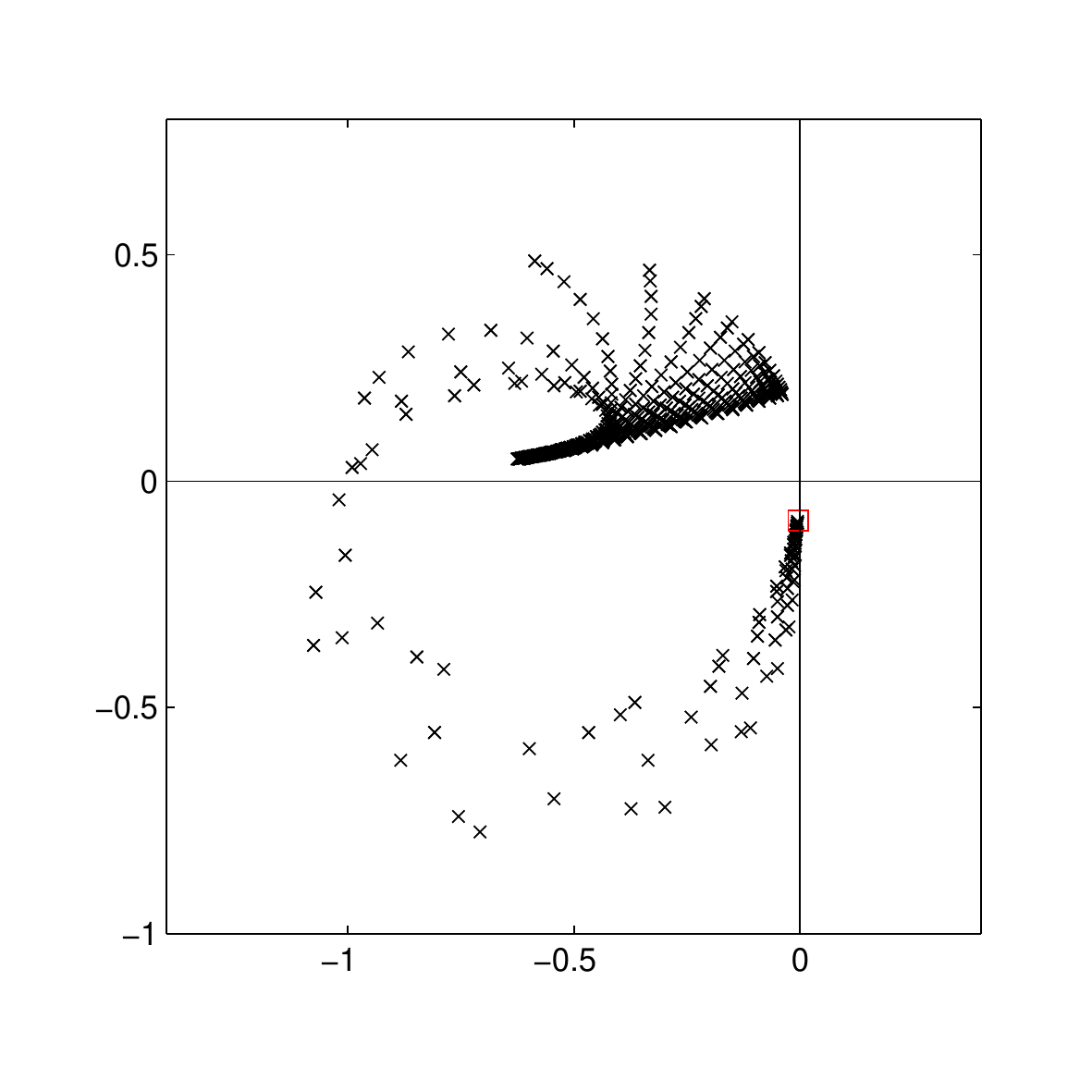} 
\includegraphics[width=6cm]{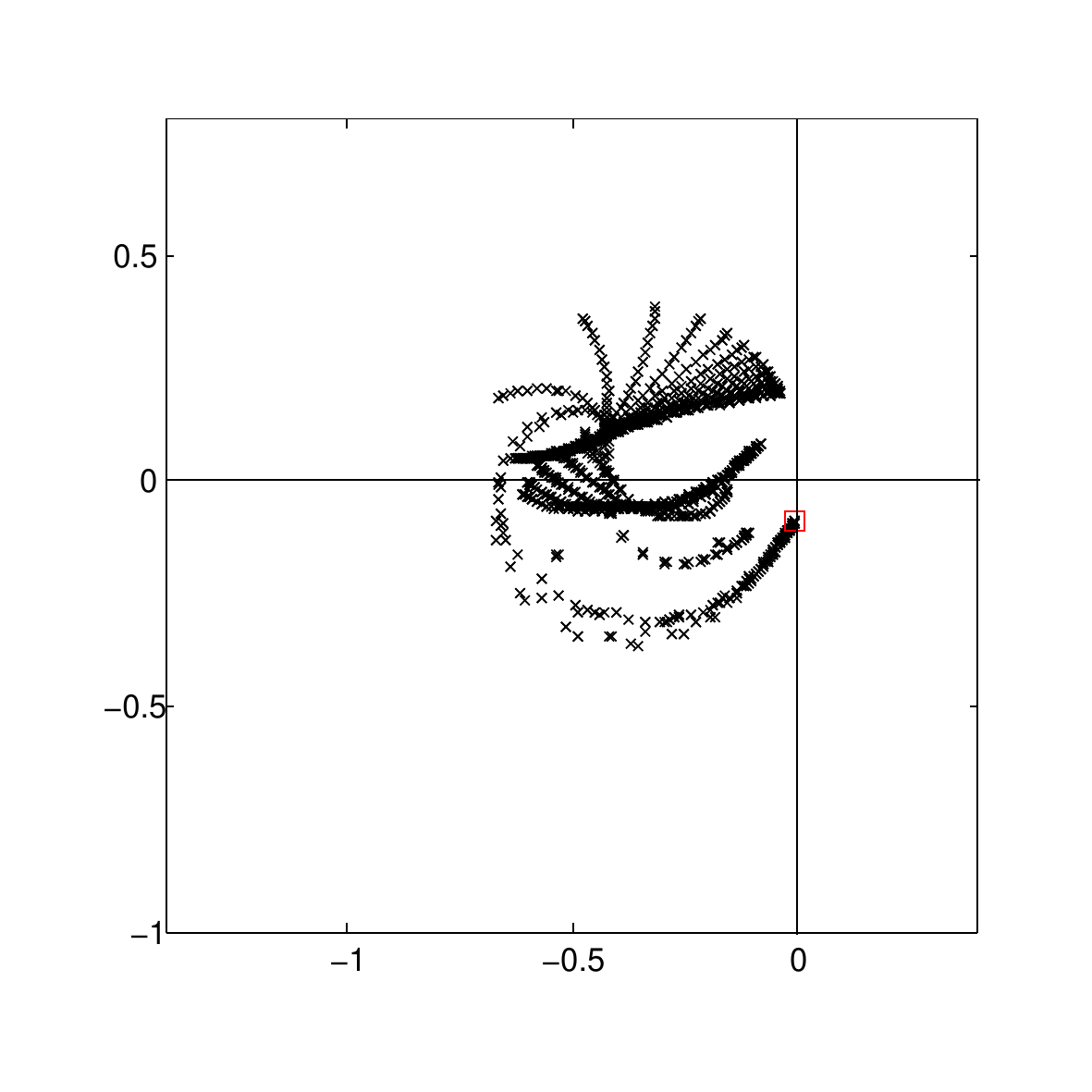} 
\caption{Spectrum of operator (\ref{eq:errchar}) displaying the eigenvalues $\gamma_l^{2h}$ (\ref{eq:gam}) for a 2D constant-$k$ Helmholtz problem with Dirichlet (left) and ECS (right) boundary conditions. Top: wavenumber $k = 20$ and $n = 32$ discretization points. Mid: wavenumber $k = 20$ and $n = 64$ discretization points. Bottom: wavenumber $k = 40$ and $n = 64$ discretization points.}
\label{fig:specgam}
\end{center}
\end{figure}

A similar spectral analysis can be performed from a slightly different point of view. In the following we focus without loss of generality on the coarse grid residual equation of the $TG^1$ correction scheme
\beq \label{eq:erreq1}
\tilde{A}^{2h} \tilde{e}^{2h} = \tilde{r}^{2h}, 
\eeq
where the right-hand side $\tilde{r}^{2h} = c \, r^{2h} = c \, I_h^{2h} r^h$ $(c \in \mathbb{C})$ is the (scaled) restricted fine grid residual $r^h = f^h - A^h v^h$. The scaling of the right-hand side is mandatory for the CSG variant of the scheme, and can intuitively be seen as a natural effect to the rotation of the spectrum, i.e.~typically $c = e^{-i d\theta}$, see also (\ref{eq:csghelmprime}). For the CSL variant, no scaling is required as the spectrum is merely translated in the complex plane, thus it suffices to set $c = 1$. Equation (\ref{eq:erreq1}) is solved for $\tilde{e}^{2h}$ which is then interpolated to correct the fine grid solution. Comparing to the standard two-grid correction scheme with residual equation
\beq 
A^{2h} e^{2h} = r^{2h}, 
\eeq
the solution $\tilde{e}^{2h}$ to the perturbed residual equation (\ref{eq:erreq1}) differs slightly from $e^{2h}$ due to the minor difference in the level-dependent scheme between the fine- and coarse grid operators. The relation between $\tilde{e}^{2h}$ and $e^{2h}$ is given by equating the above, i.e.~$\tilde{A}^{2h} \tilde{e}^{2h} = c \, A^{2h} e^{2h}$, which implies
\begin{align} \label{eq:etil}
\tilde{e}^{2h} &= c \, \tilde{A}^{{2h}^{-1}} A^{2h} \, e^{2h} = c \, \tilde{A}^{{2h}^{-1}} (A^{2h} + c^{-1} \, \tilde{A}^{2h} - c^{-1} \, \tilde{A}^{2h}) \, e^{2h} \notag \\
							 &= e^{2h} + \left(\tilde{A}^{{2h}^{-1}} (c \, A^{2h}-\tilde{A}^{2h})\right) \, e^{2h}.
\end{align}
As briefly pointed out in the introduction to this section, the level-dependent scheme indeed introduces a small additional error $\varepsilon^{2h}$ on the correction due to the difference in problem definition between coarse and fine grid. This error can effectively be characterized using equation (\ref{eq:etil}) as
\beq \label{eq:errchar}
\varepsilon^{2h} =  \left(\tilde{A}^{{2h}^{-1}} (c \, A^{2h}-\tilde{A}^{2h})\right) \, e^{2h} ,
\eeq
and preferably equals zero. In reality, however, $e^{2h}$ consists of a linear combination of coarse grid eigenmodes, which implies that the additional error $\varepsilon^{2h}$ is given by
\beq \label{eq:errlc}
e^{2h}  = \sum_{l=1}^{n/2} a_l w_l^{2h} \quad \Rightarrow \quad \varepsilon^{2h} = \sum_{l=1}^{n/2} \gamma_l^{2h} a_l w_l^{2h},
\eeq
for coefficients $a_l \in \mathbb{R}$. The weights $\gamma_l^{2h}$ are the eigenvalues of the matrix operator $\tilde{A}^{{2h}^{-1}} (c\,A^{2h}-\tilde{A}^{2h})$, given explicitly by
\beq \label{eq:gam}
\gamma_l^{2h} = \frac{ c \, \lambda_l^{A^{2h}} - \lambda_l^{\tilde{A}^{2h}} }{ \lambda_l^{\tilde{A}^{2h}} },  \qquad l = 1,\ldots,\frac{n}{2}.
\eeq
These eigenvalues can be seen as weighting each eigenmode component of the additional error $\varepsilon^{2h}$ implied by the perturbation of the problem. They are shown in Figure \ref{fig:specgam} for various problem settings. 

For the CSL level-dependent scheme without additional right-hand side scaling $(c = 1)$, the eigenvalues $\gamma_l^{2h}$ can be rewritten explicitly as
\beq \label{eq:gamk}
\gamma_l^{2h} = \frac{\lambda_l^L-k^2 - \lambda_l^L+k^2 e^{id\theta}}{\lambda_l^L -k^2e^{id\theta}} = \frac{k^2 (e^{id\theta} - 1)}{\lambda_l^L -k^2e^{id\theta}},
\eeq
where the symbol $L$ is now used to denote the discretized \emph{coarse grid} Laplacian $L = -\Delta^{2h}$. Note that for the Poisson problem with $k^2 = 0$, the above expression reduces to
\beq \label{eq:eigzero}
\gamma_l^{2h} = 0 , \qquad l = 1,\ldots,\frac{n}{2},
\eeq
hence considering (\ref{eq:errlc}) we have $\varepsilon^{2h} = 0$. This result implies that for a Poisson problem the level-dependent correction scheme is theoretically identical to the original multigrid correction scheme, as no additional error is added by solving the equation on a perturbed coarser level.

Moving back to the general Helmholtz setting where $k^2$ is distinctly different from zero, we consider a very smooth coarse grid eigenmode $w_l^{2h}$ with $l \ll n/4$ which implies $\lambda_l^L \ll k^2$. For such a smooth eigenmode it follows from (\ref{eq:gamk}) that 
\beq \label{eq:glim}
\gamma_l^{2h} \approx e^{-id\theta}-1 , \qquad l \ll \frac{n}{4}.
\eeq
Considering (\ref{eq:errlc}) and presuming that the per-level perturbation $d\theta$ is small, this result implies that smooth eigenmodes only marginally contribute to the error component $\varepsilon^{2h}$. Consequently, the additional error component $\varepsilon^{2h}$ that was induced by the adaptation of the standard correction scheme to a level-dependent scheme consists mainly of \emph{oscillatory} eigenmodes. As the elimination of the oscillatory eigenmodes is exactly the function of the pre- and postsmoothing steps encapsulating the correction scheme, the entire smoothed level-dependent two-grid scheme is expected to perform (nearly) as efficiently as the standard two-grid schemes which are currently being used as preconditioners. Moreover, the coefficients $a_l$ for oscillatory modes $w_l^{2h}$ in expression (\ref{eq:errlc}) can in fact a priori be assumed small due to presmoothing. We again stress, however, that the newly developed level-dependent scheme can be used \emph{directly} as a solver for the indefinite problem, instead of suiting preconditioning purposes. This will be demonstrated in Section 4.

We additionally note that, given a sufficiently fine discretization with respect to the wavenumber $k$ (e.g.~respecting the $kh < 0.625$ criterion, see \cite{bayliss1985accuracy}), the spectral radius of the eigenvalues $\gamma_l^{2h}$ is bounded as a function of the number of grid points $n$. Given a fixed wavenumber $k$ and perturbation angle $d\theta$, the spectral radius is maximal for the eigenvalue $\gamma_l^{2h}$ minimizing the denominator $|\lambda_l^L-k^2e^{id\theta}|$ of (\ref{eq:gamk}). This minimal distance over all $l$, however, remains quasi unchanged for increasingly fine meshes, as the spectrum of the Laplacian $L$ only expands along the positive real line as $n$ grows. As a consequence of this observation, good scalability in function of $n$ can be expected for the level-dependent scheme.

A similar analysis can be performed for the CSG level-dependent scheme, yet to obtain comparable properties to those following from expression (\ref{eq:gamk}), the right-hand side $r^{2h}$ should be scaled by $c = e^{-id\theta}$, i.e.~the same rotation used to perturb the matrix operator $\Delta^{2h}$. Note that this scaling is quite natural, see (\ref{eq:csghelmprime}), and it does not complicate the coarse grid system; however, it will make sure the Poisson problem is solved identically by the standard- and CSG level-dependent schemes, cfr.~(\ref{eq:eigzero}).

For the CSG level-dependent scheme with right-hand side scaling parameter $c = e^{-id\theta}$, the eigenvalues $\gamma_l^{2h}$ (\ref{eq:gam}) can be rewritten as
\beq \label{eq:gamk2}
\gamma_l^{2h} = \frac{\lambda_l^L e^{-id\theta} - k^2 e^{-id\theta} - \lambda_l^L e^{-id\theta} + k^2 }{\lambda_l^L e^{-id\theta}-k^2} = \frac{ k^2 (1-e^{-id\theta}) }{\lambda_l^L e^{-id\theta}-k^2},
\eeq
which is identical to expression (\ref{eq:gamk}). Consequently all properties stated above hold: the CSG level-dependent scheme with $c = e^{-id\theta}$ is identical to the standard correction scheme when solving the Poisson problem, and the additional error $\varepsilon^{2h}$ that arises when applying the level-dependent scheme to a general Helmholtz problem consists mainly of oscillatory modes, which are consecutively eliminated by the application of the smoother.

The eigenvalues $\gamma_l^{2h}$ characterizing the error $\varepsilon^{2h}$ between the original and the level-dependent multigrid schemes are plotted in Figure \ref{fig:specgam}. The smooth mode limit value $(e^{-id\theta}-1)$ from (\ref{eq:glim}) is marked by the $\square$-symbol. The eigenvalue corresponding to the smoothest eigenmode is indeed found close to zero as suggested by the discussion above. The oscillatory eigenmodes clearly have the largest contribution to the perturbation error $\varepsilon^{2h}$. Note that the spectrum does not expand as a function of the number of discretization points $n$. The numerical experiments covered in Section 4 indeed confirm that scalability in function of the number of grid points is maintained. As shown by Figure \ref{fig:specgam}, the spectrum does, however, grow with the wavenumber $k$, implying perfect $k$-scalability cannot be ensured for the level-dependent scheme.\footnote{The authors are currently unaware of any iterative method for the general Helmholtz problem which ensures perfect $k$-scalability. Significant efforts towards this aim for the class of preconditioned Krylov solvers were recently made by Vuik et al. \cite{sheikh2011scalable}, combining the CSL preconditioner with multigrid deflation. Additionally, the efforts by Brandt \& Livschitz in \cite{brandt1997wave} for constant-$k$ problems and more recent work by Engquist \& Ying \cite{engquist2011sweeping} on so-called `sweeping preconditioners' do show very promising and (near-)optimal $k$-scalability on the model problems treated therein.} The effects caused by the growing spectral radius in function of $k$ are, however, somewhat extenuated by the following two observations. First, the growing addition of oscillatory modes by the correction scheme is generally not problematic, as they are efficiently eliminated by the smoother, and second, the per-level perturbation $d\theta$ \emph{decreases} in function of the number of grid levels, such that the contribution of the smooth modes to $\varepsilon^{2h}$ (which \emph{would} be problematic if growing) diminishes in function of $n$. As a final remark we note from Figure \ref{fig:specgam} that the spectrum is much more clustered around the origin when using non-Dirichlet boundary conditions like e.g.~ECS, Sommerfeld or PML boundaries, which is particularly important regarding the contribution of smooth modes to $\varepsilon^{2h}$. Indeed, experiments have shown that the level-dependent scheme performs very poorly on a purely Dirichlet-bounded problem. As often suggested in the literature, the use of Dirichlet boundary conditions implies a `worst case scenario' for efficient solution, as no natural damping occurs (cfr.~Section 2.4). Hence, it must be noted that the level-dependent scheme is mainly devised for use on problems with more realistic ECS or Sommerfeld boundary conditions.

\subsection{On the relaxation method or `smoother'}

As the second component of the multigrid scheme, the smoother plays a crucial role in the elimination of highly oscillatory modes in the error. The smoothing of oscillatory eigenvalues is fundamental for the performance of the multigrid scheme, as it complements the action of the coarse grid correction scheme constructed above. 

\subsubsection{Introduction of polynomial smoothers.} Being a basic iterative scheme, the smoother can be chosen from a wide variety of standard relaxation methods, the most common being weighted Jacobi relaxation given by the iteration matrix
\beq \label{eq:wjacobi}
R_{\omega} = I - \omega D^{h^{-1}} A^h,
\eeq
where $\omega \in \mathbb{C}$. Note that in case of a constant diagonal $D^h$, the above expression can be rewritten as the evaluation of a first order polynomial $p_1(t)$ in $A$, where
\beq \label{eq:wjacobipol}
p_1(t) = I - \omega_1 t,
\eeq
with $\omega_1 \in \mathbb{C}$. This observation led to the development of so-called \emph{polynomial smoothers}, introduced in \cite{varga1962matrix} and first embedded in multigrid in e.g.~\cite{bank1985sharp}, in which a higher-order polynomial $p_m(t)$ of degree $m \geq 1$ is employed as a smoother
\beq
p_m(t) = (I - \omega_1 t)(I - \omega_2 t)\ldots(I - \omega_m t) = I + \sum_{i=1}^m c_i t^i,
\eeq
where $\omega_1,\ldots,\omega_m \in \mathbb{C}$ and $c_1,\ldots,c_m \in \mathbb{C}$. Note that per definition $p_m(0) = 1$. As described in \cite{reps2011multigrid}, the following two features are essential properties of any smoothing polynomial $p_m(t)$:
\begin{itemize}
\item[(1)] Stability property: $\forall h, ~ \forall \lambda_l \in spec(A^h) : |p_m(\lambda_l)| < 1 \quad(l = 1,\ldots,n)$,
\item[(2)] Smoothing property: $p_m(\lambda_n) = 0$.
\end{itemize}
The first condition expresses the smoother $p_m(t)$ is preferably stable on all levels, meaning none of the eigenmodes are amplified by its application. The second condition implies the smoother $p_m(t)$ effectively damps the oscillatory eigenmodes, mapping the corresponding high-frequency eigenvalues onto (or close to) zero.

In view of these conditions we suggest GMRES(m) as a replacement for the standard smoothing schemes, where the degree $m$ is typically chosen $m = 3$. Indeed, this implies the intrinsic construction of a smoothing polynomial of degree $m$ for which the coefficients $c_1,\ldots,c_m$ are optimally chosen with respect to a minimization of the current residual. Note that the use of GMRES(m) implies a \emph{level-dependent smoother} is included in the multigrid scheme, as the polynomial coefficients are redetermined upon every application of the smoother. 

\subsubsection{Motivation for the use of GMRES(m).} Although the application of a polynomial smoother $p_m(t)$ is computationally more expensive than standard weighted Jacobi or Gauss-Seidel relaxation, the use of an $m$-th order polynomial with variable polynomial coefficients yields significant advantages over the fixed first-order Jacobi smoothing polynomial (\ref{eq:wjacobipol}). Indeed, consider the $\omega$-Jacobi smoother with fixed relaxation parameter $\omega \in \mathbb{C}$ as given by (\ref{eq:wjacobi}). Stability of the smoother is guaranteed whenever 
\beq \label{eq:jacstab}
|1 - \omega \lambda_l^h| < 1, \qquad \forall l = 1,\ldots,n, \qquad \forall h,
\eeq
on all levels in the multigrid hierarchy. Let us, without loss of generality, assume that the problem under consideration has Dirichlet bounds, implying all eigenvalues of $A^h$ are real. Note that on the finest grid ($h$-level) the eigenvalues $\lambda_l^h$ corresponding to the oscillatory eigenmodes $w_l^h$ $(l > n/2)$ are distinctly situated on the positive real line, i.e.~$\lambda_l^h \gg 0$, implying one at least requires $\mathcal{R}(\omega) > 0$ in order to obtain a stable smoother. However, when the wavenumber $k$ is significantly larger than zero, there notably exists a coarse $H$-level in the grid hierarchy on which all eigenvalues $\lambda_l^H$ are negative, meaning $\forall l : \lambda_l^H \ll 0$. On this level one clearly requires $\mathcal{R}(\omega) < 0$ for (\ref{eq:jacstab}) to hold, resulting in a contradictory overall requirement on the relaxation parameter $\omega$. As a direct generalization of this result one can state that fixing the coefficients of a polynomial smoother inevitably results in smoother instability on certain levels in the multigrid hierarchy. This stability issue is directly resolved by using GMRES(m), which selects different coefficients on each level.

In addition to the stability issues of a polynomial smoother with fixed coefficients, one can furthermore show that for any given level in the multigrid hierarchy, a first (and even second) order polynomial can never lead to a stable smoother. Due to the indefinite nature of the problem the smoothest eigenvalues $\lambda_l^h$ $(l \ll n/2)$ on the finest levels are located on the negative real line, i.e.~$\lambda_l^h \ll 0$, whereas for the oscillatory eigenvalues $\lambda_{l'}^h$ $(l' > n/2)$ one has $\lambda_{l'}^h \gg 0$. Note that for the smoother to be stable, it follows from (\ref{eq:jacstab}) that $\forall l : \mathcal{R}(\omega) \lambda_l^h \in \,\, ]0,2[$. For positive eigenvalues $\lambda_l^h > 0$ this implies that $0 < \mathcal{R}(\omega) < 2/\lambda_l^h$, whereas for $\lambda_l^h < 0$ one requires $2/\lambda_l^h < \mathcal{R}(\omega) < 0$, which are clearly contradictory requirements. Consequently, no choice for $\omega\in \mathbb{C}$ can possibly lead to a stable first order polynomial smoother. 

One of the first papers to partially substitute the smoother by GMRES(m) is notably \cite{elman2002multigrid}. In \cite{calandra2011two} a GMRES(m) method was used as a smoother for the entire multigrid cycle. In \cite{reps2011multigrid} it was shown that for a Helmholtz problem with ECS boundaries a stable particular polynomial $p_m(t)$ of degree $m = 3$ which effectively maps the most oscillatory eigenvalues onto zero can always be constructed. Since the experiments conducted further on in this paper apply ECS boundaries to represent outgoing wave boundary conditions, we have chosen to use GMRES(3) as a smoother substitute throughout this work. Note that when using GMRES(m) as a smoother substitute within a preconditioning multigrid method, the multigrid cycle intrinsically is a variable nonlinear preconditioner which must be combined with a \emph{flexible} outer Krylov subspace method, see \cite{saad1993fgmres,simoncini2003flexible}.

\section{Numerical experiments}

In this section we extensively test the convergence and scalability of the new level-dependent multigrid scheme with respect to the number of fine grid discretization points $n$, the wavenumber $k$ and the per-level perturbation parameter $d\theta$. Additionally, we illustrate the convergence speed of the level-dependent method and compare with current state-of-the-art solution methods like Complex Shifted/Stretched preconditioned GMRES in terms of iteration count and CPU time. The three presented benchmark problems cover a wide variety of possible Helmholtz settings, providing an adequate testing framework for the new level-dependent multigrid solution method. All experiments are conducted in a 2D or 3D setting and use ECS absorbing boundary conditions, unless explicitly stated otherwise.\footnote{Hardware specifications: Intel$^{{\scriptsize\textregistered}}\hspace{-0.05cm}$ Core$^{\text{\tiny TM}}$ i7-2720QM 2.20GHz CPU, 6MB Cache, 8GB RAM. Software specifications: all numerical experiments implemented in Matlab$^{{\scriptsize\textregistered}}$.}

\subsection{The constant $k$ model}

\begin{table}[t]
\centering
\begin{tabular}{c c c c c c c c c c}
\hline \vspace{-0.35cm} \\ 
  $n_x \times n_y$       & $32^2$ & $64^2$ & $128^2$ & $256^2$ & $512^2$ & $1024^2$ \\
  	\vspace{0.1cm}
  $d\theta$ 	& $\frac{\pi}{30}$ & $\frac{\pi}{36}$ & $\frac{\pi}{42}$ & $\frac{\pi}{48}$ & $\frac{\pi}{54}$ & $\frac{\pi}{60}$ \\
\hline
  {\footnotesize MG iter}				& 77 & 33 & 25 & 25 & 25 & 28  \\
  {\footnotesize CPU total}			& 0.96 & 1.12 & 2.81 & 12.78 & 53.78 & 248.8   \\
  {\footnotesize CPU / 1000 pts} 			& 0.94 & 0.27 & 0.17 & 0.20 & 0.21 & 0.24   \\
\hline
\end{tabular}
\vspace{0.2cm}
\caption{Multigrid performance of the level-dependent scheme for a 2D constant $k$ Helmholtz problem with $k = 40$. Listed are the number of iterations, total CPU time until convergence and CPU time per 1000 grid points (in s.) for different levels of discretization. A series of V(1,1)-cycles with GMRES(3) smoothing is used as a solver.}
\vspace{-0.3cm}
\label{tab:constkh1}
\end{table}

\begin{table}[t]
\centering
\begin{tabular}{c c c c c c c c c}
\hline \vspace{-0.35cm} \\ 
  $n_x \times n_y$       & $64^2$ & $128^2$ & $256^2$ & $512^2$ & $1024^2$ & $2048^2$ \\
	  \vspace{0.1cm}
  $d\theta$ 	& $\frac{\pi}{36}$ & $\frac{\pi}{42}$ & $\frac{\pi}{48}$ & $\frac{\pi}{54}$ & $\frac{\pi}{60}$ & $\frac{\pi}{66}$ \\
\hline
  {\footnotesize MG iter}				& 180 & 57 & 39 & 40 & 40 & 43 \\
  {\footnotesize CPU total}			& 5.78 & 6.26 & 19.99 & 86.45 & 357.2 & 1552 \\
  {\footnotesize CPU / 1000 pts} 			& 1.41 & 0.38 & 0.31 & 0.33 & 0.34 & 0.37 \\
\hline
\end{tabular}
\vspace{0.2cm}
\caption{Multigrid performance of the level-dependent scheme for a 2D constant $k$ Helmholtz problem with $k = 80$. Listed are the number of iterations, total CPU time until convergence and CPU time per 1000 grid points (in s.) for different levels of discretization. A series of V(1,1)-cycles with GMRES(3) smoothing is used as a solver.}
\vspace{-0.3cm}
\label{tab:constkh2}
\end{table}

The first model problem is the most elementary Helmholtz equation and is therefore a natural subject for theoretical analysis and a benchmark for numerical experiments. It describes a general scattering problem in a homogeneous medium, i.e.\ with a constant wavenumber $k$ and a point-source located at the center of the domain $\Omega$. The solution vanishes towards infinity in all directions. We use the following standard setting of parameters:
\begin{eqnarray*}
\Omega &=& (0,1)^2,\\
k(x,y) &=& k^2, \\
\chi(x,y) &=& \begin{cases}
	1, \quad \text{for } x=y=1/2,\\
	0, \quad \text{elsewhere}, \\
\end{cases} \\
u &=& \mbox{ outgoing on }\partial\Omega.
\end{eqnarray*}
In a dimensionless formulation, the computational domain $\Omega$ is a unit square with absorbing boundary conditions on all four edges. The domain is discretized using $n$ equidistant grid points in every spatial dimension, with boundary conditions implemented by an ECS layer consisting of $n/4$ grid points surrounding the domain.

Table \ref{tab:constkh1} and \ref{tab:constkh2} show level-dependent multigrid method convergence results for the constant-$k$ model problem. The number of multigrid V(1,1)-cycles and corresponding CPU time required to solve the problem to a residual tolerance of $10^{-7}$ are given in function of the fine grid discretization. Note that a fine grid number of $n = 2^p$ discretization points corresponds to exactly $p$ levels in the multigrid hierarchy for a full V-cycle. As can be read in standard textbooks \cite{briggs2000multigrid,trottenberg2001multigrid}, scalability in function of $n$ (sometimes referred to as `$h$-scalability') is a requirement for any multigrid method, i.e.~the number of iterations is expected to remain constant in function of the fine level discretization. For the new level-dependent method one indeed observes that scalability in function of $n$ is guaranteed. The total CPU time scales with the number of fine level grid points $n$, as doubling the number of grid points implies CPU time $\times$ $2^d$.

Note how convergence deteriorates at very rough fine grid discretizations for which $kh > 0.625$, i.e.~we do not meet the requirements of 10 grid points per wavelength on the finest level of the multigrid hierarchy. Is is well-known (see \cite{bayliss1985accuracy}) that from a physical point of view traditional solution methods require a minimum number of grid points per wavelength to accurately represent (all eigenmodes of) the solution on the finest level. Additionally, this condition translates into a direct theoretical requirement for the level-dependent scheme due to the definition of the per-level perturbation parameter $d\theta = \theta_{max}/p$ where $p$ is the total number of levels in the hierarchy and $\theta_{max}$ is a fixed rotation angle (commonly chosen $\theta_{max} = \pi/6$). Consequently, as the coarse grid correction requires $d\theta$ to be sufficiently small, see (\ref{eq:glim}), convergence of the level-dependent scheme clearly benefits from a fairly fine discretization. We refer to Table \ref{tab:constkth} for a more elaborate discussion on the impact of the perturbation parameter $d\theta$ on convergence.

\begin{table}[t]
\centering
\begin{tabular}{c c c c c c c c c}
\hline
  $k$       & $20$ & $40$ & $80$  & $160$ & $320$ \\
  $n_x \times n_y$       & $32^2$ & $64^2$ & $128^2$ & $256^2$ & $512^2$ \\
\hline
  {\footnotesize MG-FGMRES} 		& 19 (0.37) & 29 (1.20) & 53 (8.09) & 106 (125) & 204 (1605) \\
  {\footnotesize MG-FGMRES(10)} & 21 (0.38) & 30 (1.13) & 62 (7.48) & 125 (72.9) & 249 (625) \\
  {\footnotesize LVL-MG-FGMRES} 		& 19 (0.25) & 30 (1.11) & 52 (7.73) & 97 (112) & 196 (1541) \\
  {\footnotesize LVL-MG-FGMRES(10)} & 20 (0.25) & 31 (1.04) & 58 (6.88) & 117 (69.4) & 220 (586) \\
  {\footnotesize LVL-MG} 	  		& 22 (0.37) & 33 (1.12) & 57 (6.26) & 111 (56.7) & 224 (488) \\
\hline
\end{tabular}
\vspace{0.2cm}
\caption{2D Constant-$k$ problem with ECS boundary conditions. Table comparing iterations and CPU time (in s.) for standard CSG-preconditioned GMRES, level-dependent MG-preconditioned GMRES and stand-alone level-dependent multigrid. A V(1,1)-cycle with GMRES(3) smoother is used for both preconditioning and level-dependent multigrid. Wavenumber-dependent discretization following the $kh < 0.625$ criterion for a minimum of 10 grid points per wavelength.}
\vspace{-0.3cm}
\label{tab:constkk}
\end{table}

\begin{table}[t]
\centering
\begin{tabular}{c c c c c c c c c}
\hline
  $k$       & $20$ & $40$ & $80$  & $160$ & $320$ \\
  $n_x \times n_y$       & $32^2$ & $64^2$ & $128^2$ & $256^2$ & $512^2$ \\
\hline
  {\footnotesize MG-FGMRES} 		& 17 (0.27) & 36 (0.88) & 73 (5.65) & 146 (56.2) & 291 (1262) \\
  {\footnotesize MG-FGMRES(10)} & 23 (0.32) & 41 (0.88) & 77 (4.33) & 164 (35.2) & 306 (331) \\
  {\footnotesize LVL-MG} 	  		& 23 (0.30) & 36 (0.73) & 64 (3.33) & 119 (23.1) & 237 (222) \\
\hline
\end{tabular}
\vspace{0.2cm}
\caption{2D Constant-$k$ problem with Sommerfeld radiating boundary conditions. Table comparing iterations and CPU time (in s.) for standard CSG-preconditioned GMRES and level-dependent multigrid. A V(1,1)-cycle with GMRES(3) smoother is used for both preconditioning and level-dependent multigrid. Wavenumber-dependent discretization following the $kh < 0.625$ criterion for a minimum of 10 grid points per wavelength.}
\vspace{-0.3cm}
\label{tab:constkksomm}
\end{table}

Tables \ref{tab:constkk} and \ref{tab:constkksomm} compare performance of the level-dependent multigrid method, applied both as a preconditioner and as a stand-alone solver, to the current state-of-the-art CSG preconditioned Krylov methods, which are (to certain extent) equivalent to the CSL-preconditioned Krylov methods as proven in \cite{reps2010indefinite}. Note that in each outer Krylov step one V(1,1)-cycle is used to approximately solve the preconditioning system, as is common practice in the literature \cite{erlangga2006novel,erlangga2006comparison,ernst2010difficult,vangijzen2007spectral}. The LVL-MG iterations counter appoints the required number of V(1,1)-cycles to solve the problem to a residual tolerance of $10^{-7}$, whereas the (LVL)-MG-Krylov iterations represent the number of outer Krylov iterations. The computational cost of the latter is thus the combined cost of the preconditioning V-cycles and the Krylov steps.

It is clear from Table \ref{tab:constkk} that the scalability of the stand-alone LVL-MG method is comparable to that of the standard MG-FGMRES solver in terms of iteration count. However, as the LVL-MG method is a directly applicable multigrid method, the entire computational cost of the outer Krylov method\footnote{We remark that the total computational cost of a single MG-FGMRES iteration is \emph{strictly greater} than the cost of a LVL-MG solver iteration, as MG-FGMRES requires both an additional matrix-vector product plus an orthogonalization procedure to be executed in comparison to the LVL-MG method. Additionally, a similar conclusion holds for the storage costs, which is due to the storage of the Krylov base vectors, yielding a total of $m+1$ times the storage cost of the fine grid solution, whereas a stand-alone multigrid scheme uses a maximum of two times this storage.} is dropped, which is clearly reflected in the CPU time. Indeed, LVL-MG significantly outperforms the non-restarted MG-FGMRES solver in terms of computational time for large systems. The occasional increase in iterations (thus in V-cycles) is negligible compared to the computational gains from directly using level-dependent multigrid as a solver. Comparing our level-dependent multigrid solver to the fast restarted MG-GMRES(10) method, one observes an effective speed-up ranging from at least 5\% up to over +30\% for large-scale problems. Summarizing, it can be stated that the stand-alone LVL-MG method is more scalable in terms of CPU time than the current state-of-the-art preconditioned Krylov methods, in particular for large scale problems with high wavenumbers.

\begin{figure}[t]
\begin{center}
\includegraphics[width=6.5cm]{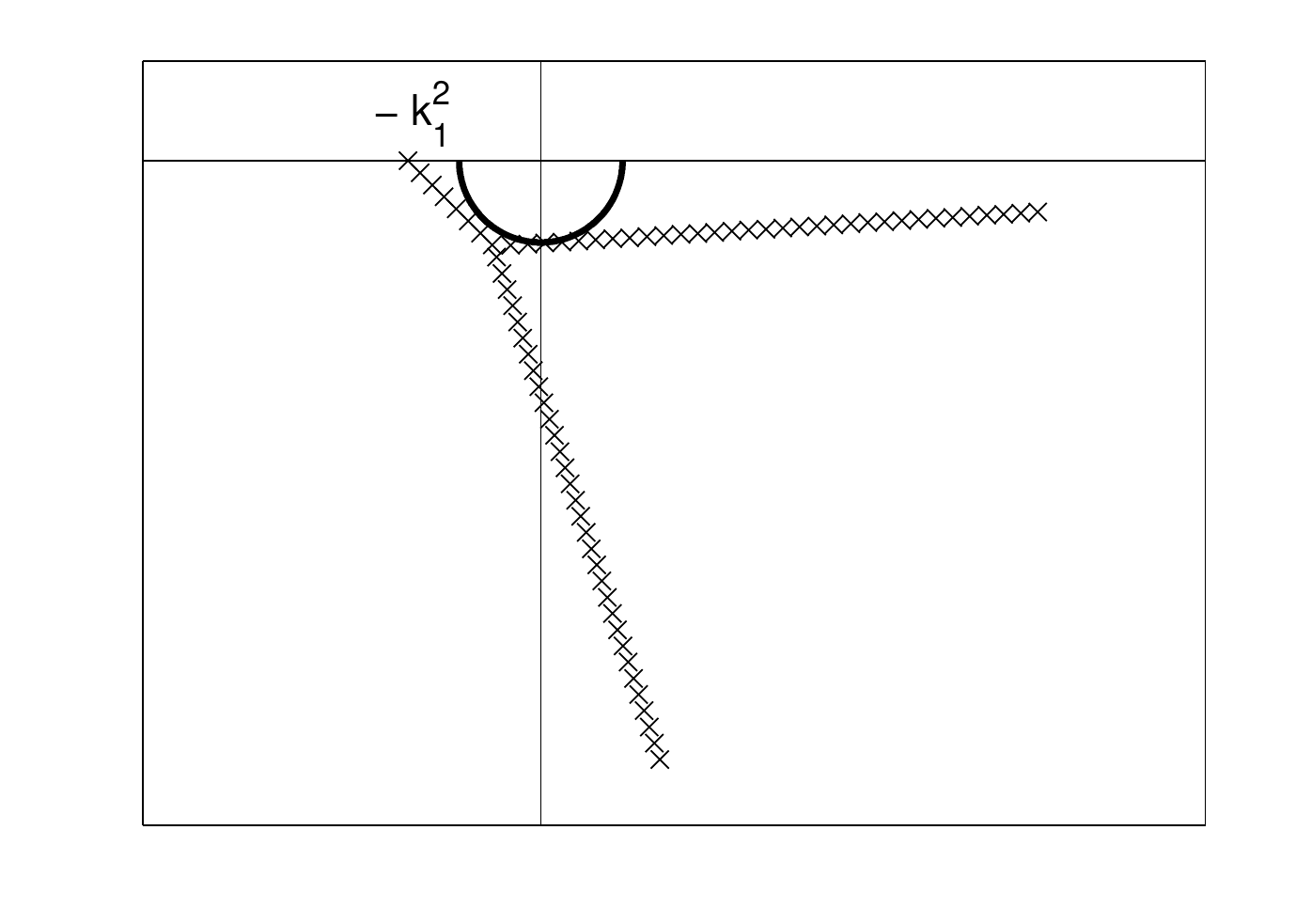}
\includegraphics[width=6.5cm]{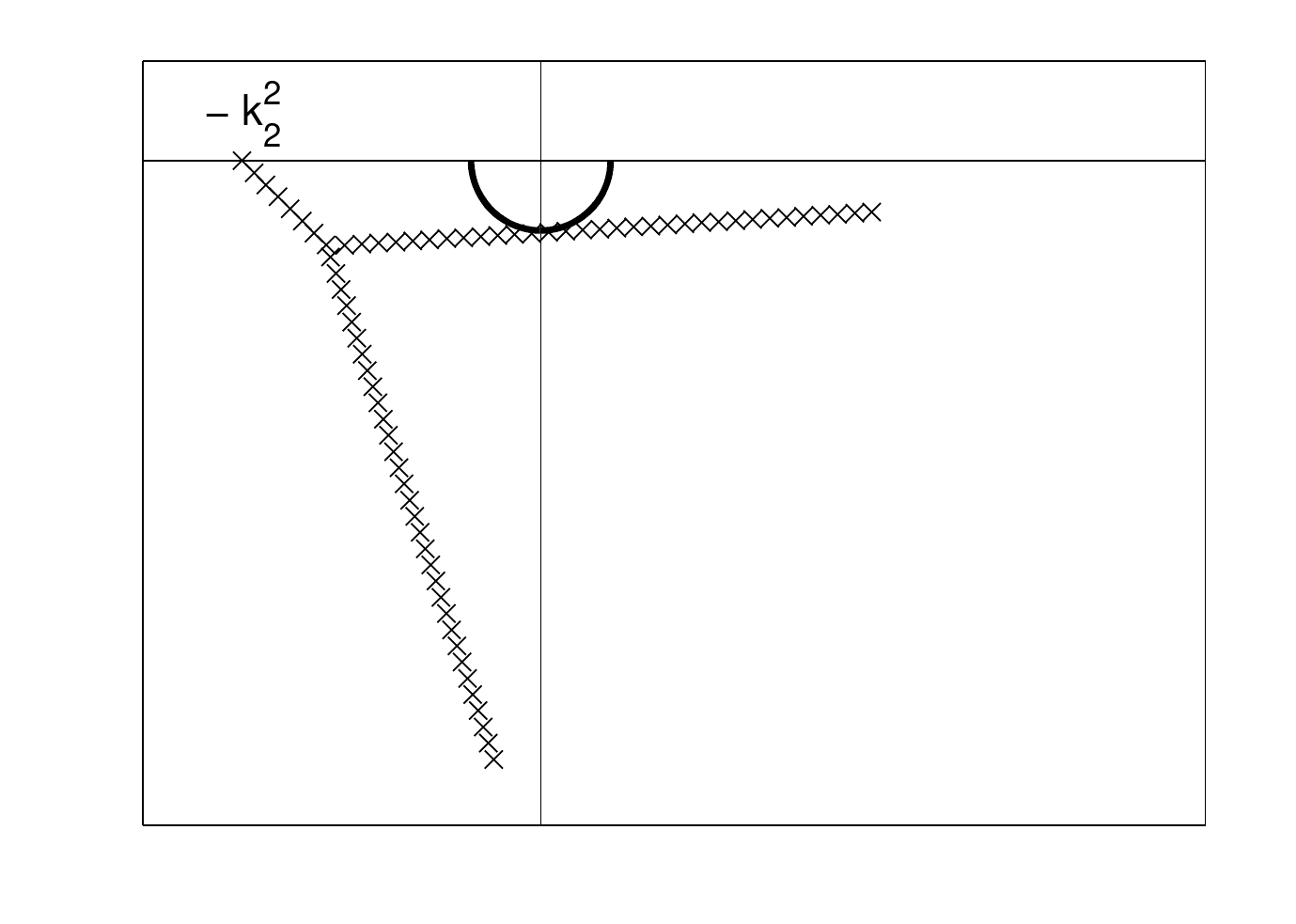}
\vspace{-0.5cm}
\caption{Schematic representation of the spectra of the perturbed ECS bounded CSG discretization matrix operator $\tilde{A}^{h}$ for different values of the wavenumber $k = k_1$ (left) and $k = k_2$ (right). Note that $k_1 < k_2$ implies that the right-hand side spectrum is generally located closer to zero compared to the left-hand side.}
\label{fig:sketchECS}
\end{center}
\end{figure}

\begin{table}[t]
\centering
\begin{tabular}{c c c c c c c c c} 
\hline
		\vspace{-0.4cm}
	\\
	  \vspace{0.1cm}
  $\theta_{max}$ 	& $\frac{\pi}{15}$ & $\frac{\pi}{12}$ & $\frac{\pi}{10}$ & $\frac{\pi}{8}$ & $\frac{\pi}{6}$ & $\frac{\pi}{5}$ & $\frac{\pi}{4}$ \\
  	\vspace{0.1cm}
  $d\theta$ 			& $\frac{\pi}{105}$ & $\frac{\pi}{84}$ & $\frac{\pi}{70}$ & $\frac{\pi}{56}$ & $\frac{\pi}{42}$ & $\frac{\pi}{35}$ & $\frac{\pi}{28}$ \\
\hline
  {\footnotesize MG iter}   & 27 	 & 25   & 23 	 & 22 	& 22 	 & 25 	& 28 \\
  {\footnotesize CPU total}	& 3.05 & 2.83 & 2.61 & 2.50 & 2.51 & 2.84 & 3.15\\
\hline
\end{tabular}
\vspace{0.2cm}
\caption{Convergence dependency of level-dependent multigrid on the per-level perturbation parameter $d\theta$. Shown are the number of iterations and CPU time (in s.) for a 2D constant $k$ Helmholtz problem with wavenumber $k = 30$ and a fine-level discretization with $n = 128$ grid points. A level-dependent V(1,1)-cycle with GMRES(3) smoother is used as a solver.}
\vspace{-0.3cm}
\label{tab:constkth}
\end{table}

When applied as a preconditioner, the performance of the new LVL-MG scheme is largely comparable to that of existing CSL/CSG-preconditioned Krylov methods. Using only a single preconditioning V(1,1)-cycle, as is common practice in multigrid preconditioning, it appears from Table \ref{tab:constkk} that applying the LVL-MG method as a preconditioner for FGMRES results in only a slight increase in performance over traditional MG-FGMRES. Note that the LVL-MG scheme was primarily designed as a solver, as it effectively solves the original problem on the finest grid (as opposed to the CSL/CSG preconditioners). The asymptotic convergence of the LVL-MG scheme to the original problem is, however, not apparent after one V-cycle, which yields a rather crude approximation to the operator inverse. Consequently, the direct application of the LVL-MG scheme as a stand-alone solver is vastly superior to its use as a preconditioner, and future numerical results will hence focus primarily on the first.

Note that for all methods good scalability in function of the wavenumber $k$ (referred to as `$k$-scalability') is de facto not guaranteed. The rising computational cost of MG-GMRES with ECS boundaries in function of $k$ was explained in \cite{reps2011analyzing} by pointing out the convergence rate behaves asymptotically when nearing $k^2 = 4/h^2$ (independently of the problem dimension). A comparable argument explains the similarly poor $k$-scalability for the level-dependent multigrid scheme. Indeed, increasing $k$ might cause the spectrum on some of the finest (only slightly perturbed) grids to contain eigenvalues which are located closer to zero, as illustrated by Figure \ref{fig:sketchECS}, resulting in a slight decline in multigrid convergence. However, the level-dependent scheme clearly has improved $k$-scalability over the MG-GMRES scheme. A combination of the new level-dependent scheme with e.g.~multigrid deflation as proposed in \cite{sheikh2011scalable} might improve $k$-scalability even further.

In Table \ref{tab:constkth} the relation between the level-dependent multigrid convergence and the perturbation parameter $d\theta$ is shown. It is clear that for very small $d\theta$ the level-dependent scheme performs poorly. Indeed, by using a too small rotation or shift some coarse level spectra might be insufficiently rotated/shifted away from zero, causing the denominator of (\ref{eq:tglvldep}) for the corresponding level(s) yet to approach zero. This possibly leads to a highly instable correction scheme, which is detrimental for the multigrid convergence. A very large per-level perturbation $d\theta$ on the other hand apparently also gives rise to worsened convergence for the level-dependent scheme due the significant difference between the fine- and coarse grid operator definition. This in particular implies that the smoothest mode weight $|1-e^{-i d\theta}|$ becomes large, see (\ref{eq:glim}), causing the additional coarse grid correction error $\varepsilon^{2h}$ to be distinctly non-zero and even contain a significant contribution of smooth modes. This corrupts the correction of the fine grid error causing convergence to deteriorate. Note that our standard choice for the maximal rotation of the coarsest level $\theta_{max} = \pi/6$ seems to be rather performant for the problem given. However, the optimal choice of the perturbation parameter $d\theta$ with respect to convergence is obviously problem-dependent.\\

\subsection{The wedge model}

\begin{figure}[t]
\begin{center}
\includegraphics[width=4cm]{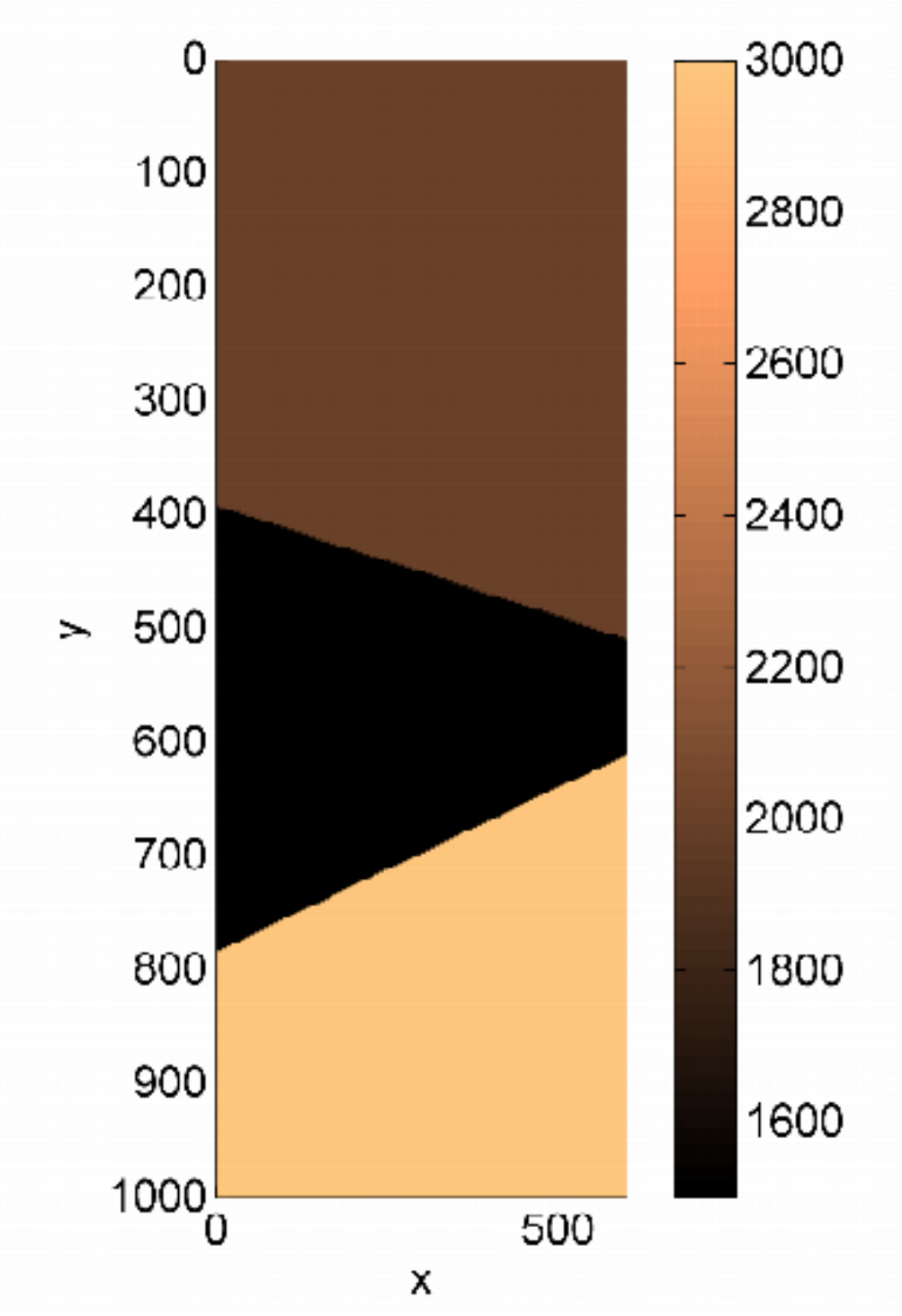}
\caption{Velocity profile $c(x,y)$ for the 2D wedge problem.}
\label{fig:wedge_velocity}
\end{center}
\end{figure}

The wedge model was introduced in \cite{plessix2003separation} for the analysis of a preconditioner based on separation of variables and adopted in \cite{erlangga2006comparison} to test the CSL preconditioner. It simulates a seismic scattering problem where a radar signal with frequency $f\in(10Hz,50Hz)$ is sent into the earth's surface that consists of three different layers. In each of these layers the sound waves travel at a different speed, as illustrated in Figure~\ref{fig:wedge_velocity},
\begin{equation*}
 c(x,y) = \begin{cases}
	2000, \quad (\text{if } 0 < y <\frac{1}{6}x+400), \\
	1500, \quad (\text{if } \frac{1}{6}x+400 \leq y < -\frac{1}{3}x+800), \\
	3000, \quad (\text{if } -\frac{1}{3}x+800 \leq y < 1000), \\
\end{cases}
\end{equation*}
which results into a mildly space-dependent wavenumber given by
\begin{equation*}
k(x,y) = \left(\frac{2\pi f}{c(x,y)}\right)^2.
\end{equation*}

\begin{table}[t]
\centering
\begin{tabular}{c c c c c c c c c}
\hline
  $f$       & $10$ & $20$ & $30$ & $40$  & $50$  \\
  $n_x \times n_y$  &  $64 \times 128$  &  $128\times 256$  &  $128\times 256$  &  $256\times 512$  &  $256\times 512$  \\
\hline
  {\footnotesize MG-FGMRES}  		& 30 (2.22) & 51 (15.3) & 78 (26.6) & 94 (218) & 114 (294) \\
  {\footnotesize MG-FGMRES(10)} & 32 (2.10) & 58 (13.8) & 85 (20.3) & 107 (128) & 133 (159) \\
  {\footnotesize LVL-MG} 		 		& 30 (2.20) & 47 (10.3) & 72 (15.2) & 83 (79.5) & 101 (96.8) \\
\hline
\end{tabular}
\vspace{0.2cm}
\caption{2D Wedge problem. Table comparing iterations and total CPU time (in s.) for standard CSG-preconditioned FGMRES and level-dependent multigrid. A V(1,1)-cycle with GMRES(3) smoother is used for both preconditioning and level-dependent multigrid. Wavenumber-dependent discretization following the $kh < 0.625$ criterion for a minimum of 10 grid points per wavelength.}
\label{tab:wedge2}
\end{table}

\noindent Additionally, the following parameters are used:
\begin{eqnarray*}
\Omega &=& (0,600)\times(0,1000),\\
\chi(x,y) &=& \begin{cases}
	1, \quad \text{for } x=300,y=0,\\
	0, \quad \text{elsewhere}, \\
\end{cases} \\
u &=& \mbox{ outgoing on }\partial\Omega.
\end{eqnarray*}

Table \ref{tab:wedge2} compares convergence results of the MG-FGMRES and LVL-MG methods on the 2D wedge problem. The observations are analogous to those made on the constant-$k$ problem: scalability with respect to the wavenumber $k$ in terms of iterations is comparable for both methods, but the new LVL-MG solver clearly outperforms the Krylov method in terms of operation count, which is reflected in the overall CPU time. Note that restarting the outer Krylov method significantly improves CPU times over regular non-restarted FGMRES due to the limited cost of the orthogonalization process. However, restarted FGMRES displays worse $k$-scalability in the number of iterations. Overall, one observes that the computational cost of the LVL-MG is considerably lower than that of the preconditioned GMRES methods (cfr.~CPU timings).

\begin{figure}[t]
\begin{center}
\includegraphics[width=4cm]{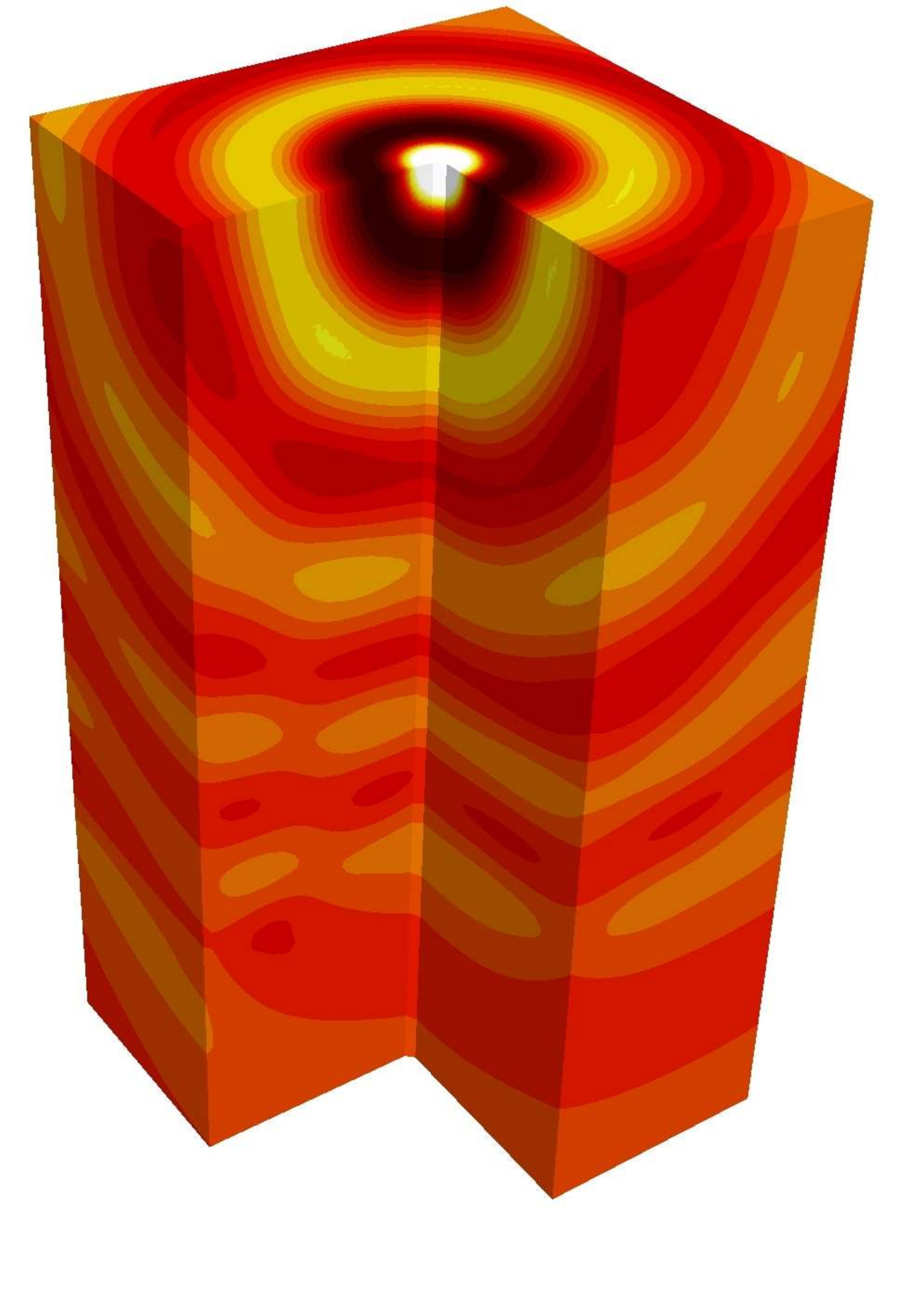} \qquad \qquad
\includegraphics[width=4cm]{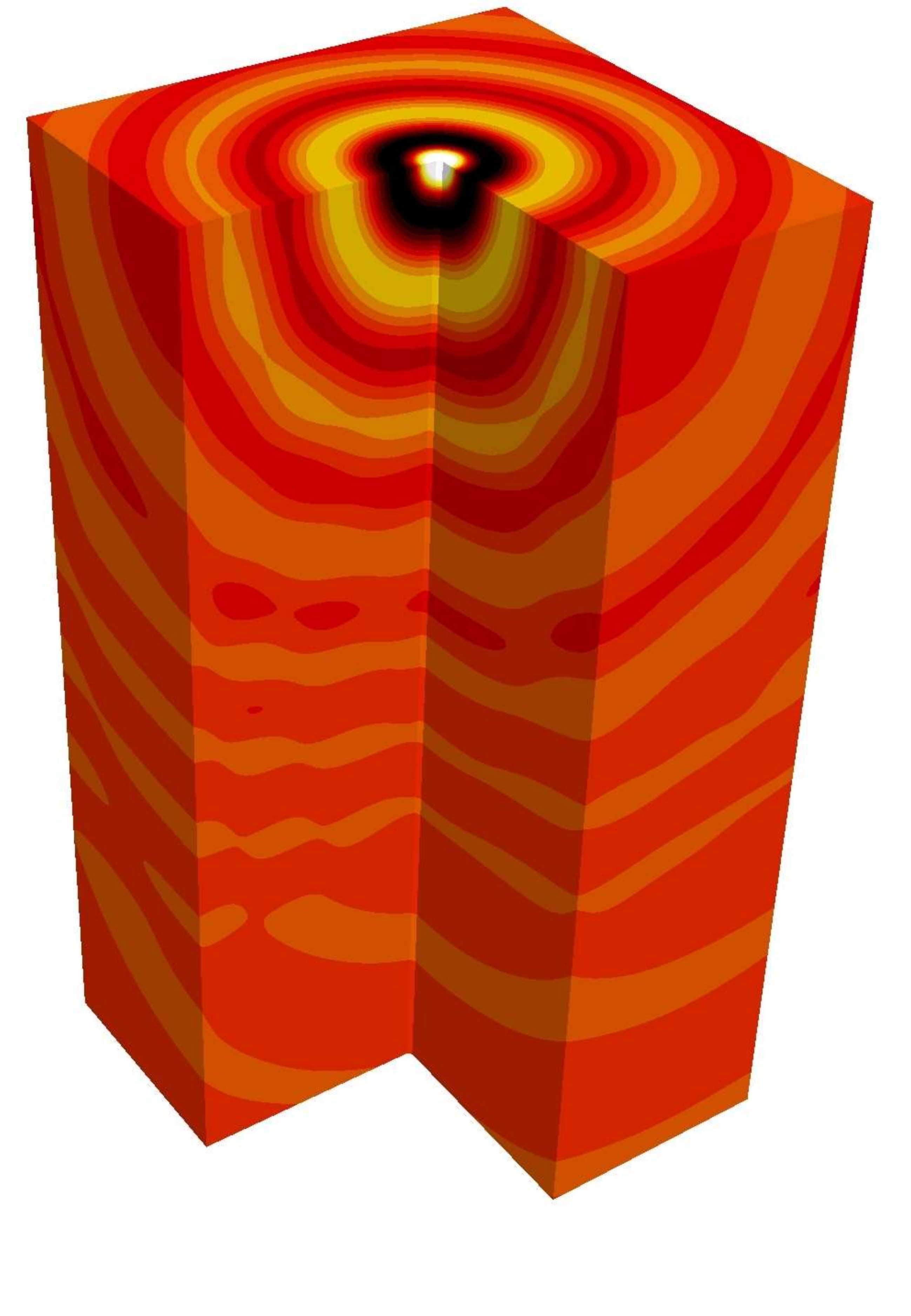} 
\caption{Real part of the solution $u(x,y,z)$ to the 3D wedge model problem for wavenumbers $k = 10$ (left) and $k = 16$ (right) respectively.}
\label{fig:wedge_solution}
\end{center}
\end{figure}

In Table \ref{tab:wedge3} the wedge problem is extended to a 3D setting. For simplicity, the third dimension (variable $z$) is chosen to be an identical copy of the $x$-dimension from the 2D formulation. Additionally, computational times have been limited by considering a moderately fine grid featuring $64\times128\times64$ grid points at solution level for all wavenumbers $k \in [10,20]$. Conclusions of the 2D problem are carried over to the 3D formulation, as iteration numbers are linearly rising in function of $k$ for both solvers. Note that the CPU time gain from using LVL-MG increases even further in higher spatial dimensions. The growing CPU time discrepancy between FGMRES and LVL-MG in function of the problem dimension is a logical consequence of the additional matrix-vector product in the outer Krylov method, which indeed gets more computationally expensive as dimension grows, and supports the use of the directly applicable level-dependent multigrid method. The resulting 3D solution $u(x,y,z)$ for wavenumbers $k = 10$ and $k = 16$ is plotted in Figure \ref{fig:wedge_solution}.

\begin{table}[t]
\centering
\begin{tabular}{c c c c c c c c c}
\hline
  $f$       & $12$ & $14$ & $16$ & $18$ & $20$  \\
  $n_x \times n_y \times n_z$ & \multicolumn{5}{c}{--------------------------~$64\times128\times64$~--------------------------} \\
\hline
  {\footnotesize MG-FGMRES}  			& 34 (294) & 39 (354) & 45 (430) & 53 (541) & 60 (645) \\
  {\footnotesize MG-FGMRES(10)} 	& 35 (234) & 40 (269) & 47 (315) & 55 (369) & 62 (416) \\
  {\footnotesize LVL-MG} 		 			& 38 (191) & 46 (232) & 50 (252) & 58 (293) & 71 (356) \\
\hline
\end{tabular}
\vspace{0.2cm}
\caption{3D Wedge problem. Table comparing iterations and total CPU time (in s.) for standard CSG-preconditioned FGMRES and level-dependent multigrid. A V(1,1)-cycle with GMRES(3) smoother is used for both preconditioning and level-dependent multigrid.}
\label{tab:wedge3}
\end{table}

\subsection{The quantum mechanical ionization model}

The Helmholtz equation can also be used to understand and predict the reaction rates of fundamental processes in few-body physics and chemistry that are of profound importance for many areas of technology. Therefore it is necessary to solve the multi-dimensional Schr\"odinger equation, equivalent to a multi-dimensional Helmholtz equation 
\beq
	-\Delta u(\bold{x})- k(\bold{x}) u(\bold{x}) = \chi(\bold{x}), \qquad \bold{x} \in \Omega \subset \mathbb{R}^n,
\eeq
with outgoing waves boundary conditions, source function $\chi(\bold{x})$ and a space-dependent wavenumber $k(\bold{x})$. 

The below 2D benchmark originates from the dynamics of two electrons in a Hydrogen molecule. Coordinates $x$ and $y$ should be interpreted as radial distances of particles to the center of mass of the system. Along the lines $x=0$ and $y=0$ the wave function $u(x,y)$ is zero, while at the other sides of the domain outgoing wave conditions are needed. We use the following model specifications:
\begin{eqnarray*}
	\Omega &=& (0,r)^2,\\
	k(x,y) &=& \frac{1}{e^{x^2}} + \frac{1}{e^{y^2}} + k_0^2, \\
	\chi(x,y) &=& \frac{1}{e^{x^2+y^2}}, \\
	u(0,y)=u(x,0) &=& 0, \qquad \quad ~ (\text{Dirichlet conditions}),\\
	u(r,y\neq0),u(x\neq0,r) &=& \text{ECS}, \qquad \text{(outgoing wave conditions)},
\end{eqnarray*}
with $0<k_0<5$. The domain $\Omega$ is chosen to range from $0$ to $r = 50$ in each spatial dimension. 

\begin{table}[t]
\centering
\begin{tabular}{c c c c c c c c c}
\hline
  $k_0$       & 1 & 2 & 3 & 4  & 5  \\
  $n_x \times n_y$  &  $128^2$  &  $256^2$  &  $256^2$  &  $512^2$  &  $512^2$  \\
\hline
  {\footnotesize MG-FGMRES}  			& 51 (5.40) & 92 (46.8) & 191 (137) & 174 (813) & 306 (2185) \\
  {\footnotesize MG-FGMRES(10)}  	& 66 (5.52) & 140 (46.7) & 245 (81.3) & 250 (421) & 393 (661) \\ 
  {\footnotesize MG-FGMRES(30)}  	& 52 (4.86) & 93 (34.7) & 226 (83.2) & 278 (565) & 426 (867) \\ 
  {\footnotesize LVL-MG} 		 			& 44 (3.42) & 83 (25.3) & 208 (63.5) & 149 (215) & 289 (418) \\
\hline
\end{tabular}
\vspace{0.2cm}
\caption{2D Ionization problem. Table comparing iterations and total CPU time (in s.) for standard CSG-preconditioned FGMRES and level-dependent multigrid. A V(1,1)-cycle with GMRES(3) smoother is used for both preconditioning and level-dependent multigrid. wavenumber-dependent discretization following the $kh < 0.625$ criterion for a minimum of 10 grid points per wavelength.}
\label{tab:ion}
\end{table}

Serving as our final test case, quantum mechanical ionization problems are generally considered rather hard-to-solve Helmholtz problems due to the heterogeneity in the domain. The wavenumber function $k(x,y)$ is clearly heavily space-dependent, rendering these type of Helmholtz problems one of the major challenges in the development of efficient and robust solvers.

Convergence results on the solution of the ionization problem using both MG-FGMRES and the LVL-MG solver are displayed in Table \ref{tab:ion}. The $k$-scalability is almost identical for both methods in terms of iteration count, as was observed in previous experiments. However, due to the elimination of the outer Krylov method when employing the level-dependent multigrid scheme, a speed-up of at least 30\% can be achieved by direct application of the LVL-MG method. Indeed, the LVL-MG method has a drastically reduced flop count in comparison to the MG-FGMRES solver, resulting in a significantly improved CPU-time scalability in function of the wavenumber $k_0$. Note that due to the use of the flexible variant of GMRES, for some model problems a small restart value $(m = 10)$ is preferable over a larger one $(m = 30)$ in terms of iterations. 

\section{Conclusions}

In this paper, we have developed a novel level-dependent correction scheme for indefinite Helmholtz problems. The proposed scheme is based on the idea of perturbing the Helmholtz operator, which was originally introduced by the Complex Shifted Laplacian and Complex Stretched Grid preconditioning techniques. These schemes respectively shift or rotate the Helmholtz operator spectrum in the complex plane away from the origin, hence guaranteeing multigrid stability. The level-dependent correction scheme incorporates this idea by gradually shifting/rotating the spectrum throughout the multigrid hierarchy, thus resolving the classical two-grid instability issues of the standard two-grid scheme while maintaining the original Helmholtz operator on the finest grid. This results in a new multigrid scheme which, contrarily to CSL or CSG, is capable of \emph{directly} solving the Helmholtz system, instead of being used as a preconditioner. 

The aforementioned methodology of gradually perturbing the original Helm\-holtz operator throughout the hierarchy introduces an unwanted `perturbation error' which is added to the coarse grid correction. This additional error is intrinsic to the level-dependent solver, as it appears as an artifact of the difference between the coarse and fine grid residual equations. However, it is shown that this perturbation error consist primarily of oscillatory eigenmodes, which are subsequently damped by the action of the smoother. Consequently, given a sufficiently performant smoother, the level-dependent two-grid scheme is expected to be a very efficient and stable solver for Helmholtz problems. 

To ensure effective smoothing we propose the use of GMRES(3) as a smoother substitute, which is recently shown to perform efficiently as a relaxation scheme in a multigrid setting. Numerical experiments on 2D and 3D Helmholtz benchmark problems of various difficulty confirm the efficiency of the proposed level-dependent correction scheme as a solver, showing it to be competitive with the current state-of-the-art CSL- or CSG-preconditioned Krylov methods. Additionally, like any multigrid solver, the level-dependent scheme is shown to be fully $h$-independent. Perfect $k$-scalability is generally not guaranteed; indeed, scalability in function of the wavenumber is largely comparable to the present-day preconditioned Krylov methods. A hybrid combination of the new level-dependent scheme presented in this paper with e.g.~multigrid deflation would be likely to address this wavenumber scalability issue.

\section{Acknowledgements}
This research has been funded by the \textit{Fonds voor Wetenschappelijk Onderzoek (FWO)} project G.0.120.08 and \textit{Krediet aan navorser} project number 1.5.145.10. Additionally, this work is partly funded by Intel$^{{\scriptsize\textregistered}}\hspace{-0.05cm}$ and by the \textit{Institute for the Promotion of Innovation through Science and Technology in Flanders (IWT)}. The authors would like to thank Hisham bin Zubair for sharing a multigrid implementation.

\nocite{*}
{\small
\bibliographystyle{plain}
\bibliography{refs}
}

\end{document}